\documentclass{article}

\usepackage{arxiv}

\usepackage[utf8]{inputenc} 
\usepackage[T1]{fontenc}    
\usepackage{hyperref}       
\usepackage{url}            
\usepackage{booktabs}       
\usepackage{amsfonts}       
\usepackage{nicefrac}       
\usepackage{microtype}      
\usepackage{lipsum}		
\usepackage{graphicx}
\usepackage{doi}

\usepackage{amsmath,amssymb,enumerate,graphicx,psfrag}
\usepackage{latexsym}
\usepackage{amsfonts}
\usepackage{curves}
\usepackage{xcolor}
\usepackage{colortbl}
\usepackage{makecell}

\usepackage{mathrsfs}
\usepackage{amsthm, amscd}
\usepackage{amsfonts,amsmath,amssymb,amsthm,cite,enumerate,epsfig,color,xcolor}
\usepackage{kotex}
\usepackage{algorithm}
\usepackage{algpseudocode}

\theoremstyle{definition}

\title{Enhanced physics-informed neural networks with domain scaling and residual correction methods
for multi-frequency elliptic problems}

\date{}

\author{Deok-Kyu Jang \\
		Department of Applied Mathematics\\ 
		Kyung Hee University\\ 
		Yongin, Republic of Korea \\
	\texttt{dkjang@khu.ac.kr} \\
	\And
	Hyea Hyun Kim \\
		Department of Applied Mathematics\\ 
		Kyung Hee University\\ 
		Yongin, Republic of Korea \\
	\texttt{hhkim@khu.ac.kr} \\
	\And
	Kyungsoo Kim$^*$ \\
		Department of Applied Mathematics\\ 
		Kyung Hee University\\ 
		Yongin, Republic of Korea \\
	\texttt{kyungsoo@khu.ac.kr} \\
}

\renewcommand{\shorttitle}{EPINN with domain scaling and residual correction methods}

\hypersetup{
pdftitle={Enhanced physics-informed neural networks with domain scaling and residual correction methodsfor multi-frequency elliptic problems},
pdfsubject={physics.CP, coRR.ML},
pdfauthor={Deok-Kyu Jang, Hyea Hyun Kim, Kyungsoo Kim},
pdfkeywords={neural network approximation, physics-informed neural networks, multi-frequency solution, domain scaling, residual correction},
}

\begin{document}
\maketitle

\begin{abstract}
In this paper, neural network approximation methods are developed for elliptic partial differential equations with multi-frequency solutions. Neural network work approximation methods have advantages over classical approaches in that they can be applied without much concerns on the form of the differential equations or the shape or dimension of the problem domain. When applied to problems with multi-frequency solutions, the performance and accuracy of neural network approximation methods are strongly affected by the contrast of the high- and low-frequency parts in the solutions. To address this issue, domain scaling and residual correction methods are proposed. The efficiency and accuracy of the proposed methods are demonstrated for multi-frequency model problems.
\end{abstract}


\section{Introduction}

A physics-informed neural network (PINN) combines the constraint-satisfaction ability of partial differential equations (PDEs) with the representation power of deep neural networks to learn solutions of PDEs.
PINNs were first introduced in~ \cite{lagaris1998artificial,mathews2021uncovering,raissi2019physics}
as a way of solving problems in mathematical physics and engineering that can be modeled as PDEs.
The idea behind PINNs is to treat the solution of a PDE as an unknown function that can be represented by a neural network.
The neural network is then trained end-to-end to satisfy the boundary conditions and PDE constraints.
This enables PINNs to deal with problems that are challenging to solve using conventional numerical techniques,
such as, those with high-dimensional input spaces and complex boundary conditions.
Due to the growing need for effective solutions to challenging physical problems in fields like fluid dynamics, structural mechanics, and heat transfer, PINNs have become increasingly popular in recent years.
Computational and theoretical studies on PINNs have also shown to be useful for problems in machine learning, computer vision, and other fields outside physics and engineering due to their flexibility and representational power.

PINNs have been applied to a variety of problems in physics, engineering, and other fields, including solving PDEs, modeling physical systems, and carrying out data-driven simulations. However, there are still some obstacles that arise when applying them to the field of computational science and engineering.
One of these obstacles is that deep neural networks (DNNs) typically perform better with data that have low-frequency contents; this is known as the Frequency Principle (F-Principle)~\cite{luo2021frequency, xu2019frequency}.
This means that DNNs can quickly learn and generalize low-frequency data, but tend to struggle when handling data with high-frequency contents~\cite{rahaman2019spectral,xu2019frequency,xu2019training}. This training behavior observed with PINNs is unlike that of the widely used multi-grid method (MGM) for solving PDEs. The MGM uses smoothing operations to capture the high-frequency spectrum of the solution first
and produce smooth and global coarse residual errors
that can be corrected efficiently by a direct solver~ \cite{mccormick1987multigrid}.

In this study, we develop a PINN approach with procedures similar to the MGM.
For the multi-frequency model solution,
we first train a neural network that can approximate its high-frequency parts effectively
and then train the remaining low-frequency parts by using a residual error equation.
To approximate the high-frequency parts effectively, we propose a novel domain scaling approach,
in which the original problem domain is scaled by a large enough factor $b>1$ to make
the high-frequency parts smooth and less oscillatory in the scaled domain. In addition,
the corresponding forcing term for the low-frequency parts, e.g., i.e., $-\triangle u$ in the Poisson problem, becomes negligibly small
in the scaled domain problem.
After forming such a scaled domain problem, we can train the neural network effectively
to learn the high-frequency parts of the original problem solution.
The trained neural network is then scaled back to the original domain to provide a neural network approximate solution to the original model problem.
In the scaled-back neural network solution, the low-frequency parts of the original model solution
are not well-resolved and they remain as dominant errors.
To resolve these remaining errors, we form a residual error equation and train
an additional neural network for the residual error equation.
By adding the scaled-back neural network solution and the neural network solution for the residual error equation,
we can form the resulting well-resolved neural network approximation to the multi-frequency model problem.

We note that some successful approaches have been developed 
to address similar problems with neural network approximations to multi-frequency models.
In the phase-shifting deep neural network (PhaseDNN)~\cite{cai2020phase},
high-frequency components of data are shifted to a low-frequency spectrum that can be effectively learned by a neural network. The original high-frequency data can then be approximately reconstructed by shifting the learned data back to their original high-frequency range. The PhaseDNN has been shown to be effective in handling highly oscillatory data from solutions of high-frequency Helmholtz equations.
However, the PhaseDNN approach results in the creation of multiple small neural networks, requiring a high computational cost for numerous phase shifts applied independently in each coordinate direction.
In the multi-scale DNN (MscaleDNN)~\cite{li2020multi,liu2020multi,wang2020multi}, high frequencies are converted into lower ones by using a radial partition of the Fourier space and a scaling-down operation.
Furthermore, it incorporates the multi-resolution concept from wavelet approximation theory, using compact scaling and wavelet functions, and replaces traditional global activation functions with ones that have compact support.
However, it still requires a relatively large network to handle the different frequencies.

In our study, we only employ two basic feed-forward neural networks with relatively small parameters
to approximate a multi-frequency model effectively, in contrast to the previous studies.
Another key ingredient in our method is the use of the $\sin(x)$ activation function; $\tanh(x)$ is the most commonly used activation function in PINNs.
The $\sin(x)$ activation function has advantages over the others, in that
it can preserve the derivative information after differentiation and is thus more suitable to approximate
PDE solutions and multi-frequency model solutions effectively; see~ \cite{sitzmann2020implicit}.
In our numerical results, the $\sin(x)$ activation function 
also provides accurate and effective neural network
approximate solutions for multi-frequency model problems.

The remainder of this paper is organized as follows.
In Section~\ref{sec2}, we include a brief introduction to PINNs and in Section~\ref{sec3}, we describe
our proposed methods, the domain scaling approach and the residual correction approach.
In Section~\ref{sec4}, numerical results are presented to demostrate the effectiveness
of the proposed methods for the multi-frequency model problems.

\section{Conventional approaches to physics-informed neural networks}\label{sec2}

In a PINN, a neural network function $\mathcal{N}(\mathbf{x};\theta)$ is used to approximate the solution of a model problem,
which often appears as a differential equation in a bounded domain $\Omega$ with a prescribed boundary condition.
For example, the solution $u(\mathbf{x})$ to the Poisson problem in a domain $\Omega$,
\begin{equation}\label{model}
\begin{split}
-\triangle u(\mathbf{x}) &=f(\mathbf{x}), \quad \mathbf{x} \in \Omega,\\
u(\mathbf{x})&=g(\mathbf{x}), \quad \mathbf{x} \in \partial \Omega,
\end{split}
\end{equation}
can be approximated by a neural network function $\mathcal{N}(\mathbf{x};\theta)$ with the parameter $\theta$ in the neural network function trained to minimize the loss value $L(\theta)$.
To give a good approximate solution, the loss value $L(\theta)$ is formed by using the residual values in the differential equation and the boundary condition,
$$L(\theta):=\frac{1}{|\Omega|}\int_{\Omega} |f(\mathbf{x}) +\triangle \mathcal{N}(\mathbf{x};\theta)|^2 \, d\mathbf{x} + \frac{1}{|\partial \Omega|}\int_{\partial \Omega} |g(\mathbf{x})-\mathcal{N}(\mathbf{x};\theta)|^2 \, ds,$$
where $|A|$ denotes the area or length of the domain $A$.
In general, finding the optimal parameter $\theta$ for the loss value $L(\theta)$ is not tractable and the following practical loss function is thus often used:
\begin{equation}\label{loss}
L_m(\theta):=w_1 \left( \frac{1}{|X_{\Omega}|} \sum_{ \mathbf{x} \in X_{\Omega} } |f(\mathbf{x})+\triangle \mathcal{N}(\mathbf{x};\theta)|^2 \right) + w_2 \left( \frac{1}{|X_{\partial \Omega}|} \sum_{ \mathbf{x} \in X_{\partial \Omega}} |g(\mathbf{x})-\mathcal{N}(\mathbf{x};\theta)|^2 \right),
\end{equation}
where $X_{A}$ denotes a set of points selected from the domain $A$ and $|X_{A}|$ is the number of points in the set $X_{A}$. In the above, $w_1$ and $w_2$ are hyperparameters to be chosen to enhance the parameter training efficiency.
We used the subscript $m$ to stress that the approximated loss value is obtained from the Monte Carlo method by using a sufficiently large number of points.

To find the optimal parameter $\theta$ to minimize the loss value $L_m(\theta)$, gradient-based methods are often used. In such approaches, the parameter training performance is strongly affected by the neural network function structure,
the training data sets $X_{\Omega}$ and $X_{\partial \Omega}$, and the hyperparameters in the loss function $L_m(\theta)$. In addition, when the model problem solution has multi-frequency features with low- and high-frequency parts, the neural network function resolves the low-frequency features at the early training stage and needs more training epochs to resolve the high-frequency part of the solution~\cite{rahaman2019spectral,xu2019frequency,xu2019training}.
Even after a large number of training epochs, the neural network approximate solution may fail to resolve the high-frequency part.

\section{Enhanced physics-informed neural networks for multi-frequency model problems}\label{sec3}
Our goal is to enhance the training efficiency and solution accuracy in neural network approximation to multi-frequency model problems. Our approach consists of two key ingredients, the domain scaling method and the residual correction method.
The effect of the domain scaling method is similar to the smoothing operator in the MGM, as it helps capture the high-frequency components of the multi-frequency model solution.
The residual correction method plays the same role as in the MGM, where the remaining low-frequency solution components
are captured by solving the residual error equation after the smoothing step.
In this section, we describe the two steps in detail.

\subsection{The domain scaling method}
We consider the model problem \eqref{model} and assume that its solution has multi-frequency features.
When the solution has high-frequency components,
the gradient-based method requires enormous training epochs to capture the high-frequency solution parts,
since some of the parameters $\theta$ in the neural network solution must have larger values to resolve the high-frequency solution parts with sufficient accuracy.
To efficiently train the high-frequency solution parts in the neural network approximation,
we propose the domain scaling method.
For a given scaling factor $b>1$, we scale the problem domain by applying a change of variables, $\widehat{\mathbf{x}}=b\mathbf{x}$, and reformulate the model problem \eqref{model}
into the following problem in the scaled domain $\Omega_b = \{b\mathbf{x} \left| \mathbf{x} \right. \in \Omega \}$:
\begin{equation}\label{scale:model}
\begin{split}
-\triangle \widehat{u}(\widehat{\mathbf{x}}) &= \widehat{f}(\widehat{\mathbf{x}}), \quad \widehat{\mathbf{x}} \in \Omega_b,\\
\widehat{u}(\widehat{\mathbf{x}})&=\widehat{g}(\widehat{\mathbf{x}}), \quad
\widehat{\mathbf{x}} \in \partial \Omega_b,
\end{split}
\end{equation}
where
$$\widehat{u}(\widehat{\mathbf{x}})=u \left( \frac{1}{b}\widehat{\mathbf{x}} \right), \quad \widehat{f}(\widehat{\mathbf{x}})= \frac{1}{b^2} 
f \left( \frac{1}{b}\widehat{\mathbf{x}} \right), \quad \text{and} \quad \widehat{g}(\widehat{\mathbf{x}})=g\left( \frac{1}{b}\widehat{\mathbf{x}} \right).$$
For the above scaled model problem \eqref{scale:model}, we approximate the solution $\widehat{u}(\widehat{\mathbf{x}})$
by a neural network function $\widehat{\mathcal{N}}(\widehat{\mathbf{x}};\widehat{\theta})$.
After training the parameters $\widehat{\theta}$ in the approximate solution $\widehat{\mathcal{N}}(\widehat{\mathbf{x}};\widehat{\theta})$,
we scale the approximate solution back to the original domain to obtain $$\mathcal{N}(\mathbf{x};\theta_s):=\widehat{\mathcal{N}}(b\mathbf{x};\widehat{\theta}),$$
where $\theta_s$ is related to the scaling factor $b$ and the trained parameters $\widehat{\theta}$.
We use the subscript $s$ to signify that it is obtained using the scaling approach.
By using the same structured neural network function $\mathcal{N}(\mathbf{x};\theta)$, we can also train the parameters $\theta$ to the model problem \eqref{model} directly and obtain the trained parameter $\theta_o$.
Here, we use the subscript $o$ to indicate that the network parameters are trained for the original model problem.

We can give a heuristic explanation of how the proposed scaled problem can yield more accurate results than the original problem when approximating the high-frequency parts of the model solution.
When the solution has high-frequency parts, some of the trained parameters $\theta_o$ in $\mathcal{N}(\mathbf{x};\theta_o)$ must take large enough values to resolve them.
In such a case, the gradient-based parameter training method
often requires enormous training epochs to produce such large values in the trained parameters $\theta_o$, and even then it can fail
to resolve parts with particularly high frequencies.
On the other hand, in the scaling approach, some of the trained parameters $\widehat{\theta}$ are magnified by the factor $b$ and the resulting scaled-back neural network solution $\mathcal{N}(\mathbf{x};\theta_s)$ retains
the advantages provided by the scaled-back large parameter values without the expense of a large number of training epochs.
Moreover, large values of the right-hand side function $f(\mathbf{x})$ in the original problem~\eqref{model} associated with the high-frequency solution parts are greatly reduced by the factor $1/b^2$ in the scaled problem~\eqref{scale:model}, and the effect of low-frequency parts diminishes in the scaled problem.
As a result, the high-frequency parts of the solution to the original problem
are smoothed out and remain as the dominant parts of the solution to the scaled problem,
and can thus be trained efficiently by neural network approximation in the scaled problem.

We now present the hyperparameter settings in the loss function to the scaled model problem~\eqref{scale:model}.
We define the loss function for this problem as
\begin{equation}\label{scale:loss:1}
L_{m,b}(\widehat{\theta}) =w_1 \left( \frac{1}{|X_{\Omega_b}|} \sum_{ \widehat{\mathbf{x}} \in X_{\Omega_b} } \left| \widehat{f}(\widehat{\mathbf{x}})+\triangle \mathcal{N}(\widehat{\mathbf{x}};\widehat{\theta}) \right|^2 \right) + w_2\left( \frac{1}{|X_{\partial \Omega_b}|} \sum_{ \widehat{\mathbf{x}} \in X_{\partial \Omega_b}} \left| \widehat{g}(\widehat{\mathbf{x}})-\mathcal{N}(\widehat{\mathbf{x}};\widehat{\theta}) \right|^2 \right).
\end{equation}
With the numbers $|X_{\Omega_b}|$ and $|X_{\partial \Omega_b}|$ kept the same for any scaling factor $b$,
the first term in the above loss function depends
on the scaling factor $b$, since the right-hand side function $\widehat{f}(\widehat{\mathbf{x}})$ is given by
$\widehat{f}(\widehat{\mathbf{x}})= 
f \left( \widehat{\mathbf{x}}/b \right) /b^2.$
When the distribution of $\widehat{f}(\widehat{\mathbf{x}})$ contains values that are too large or too small, the parameter training for the neural network approximation becomes inefficient and often fails to find a good approximate solution.
In addition, the imbalance between the differential equation residual loss value (the first term in the right-hand side of \eqref{scale:loss:1}) and the boundary condition loss value (the second term in the right-hand side of \eqref{scale:loss:1}) may prevent a well-approximated neural network solution from being found.
As a remedy, we set the weight factor
\begin{equation}\label{weight:w1}
w_1=\dfrac{1}{\Vert \widehat{f} \Vert_{X_{\Omega_b}}^2}
\quad  \text{with} \quad 
\Vert \widehat{f} \Vert_{X_{\Omega_b}}^2 = \frac{1}{\vert X_{\Omega_b} \vert} \sum_{\widehat{\mathbf{x}} \in X_{\Omega_b}} \widehat{f}(\widehat{\mathbf{x}})^2.
\end{equation}
Since the first term in the loss function is already normalized for the data distribution, we simply set
the weight factor for the boundary condition loss term to the value 1, i.e., $w_2=1$.
In our numerical computation, we observed that this simple choice for $w_2$ gives
robust and efficient training results once the first loss term is normalized using the value $\|\widehat{f}\|_{X_{\Omega_b}}^2$.

Before we move to the residual correction step, we introduce some notation that we will use later
in our numerical experiments to estimate the training and generalization performance
of the scaled neural network approximation $\mathcal{N}(\mathbf{x};\theta_s)$.
For a given data set $Y$ of the original domain $\Omega$, we define the relative error of the scaled-back neural network solution $\mathcal{N}(\mathbf{x};\theta_s)$ and its relative residual error in the $l^2$-norm as
\begin{equation}\label{rel_L2_err}
\epsilon_u (Y)= \sqrt{\frac{\sum_{\mathbf{x} \in Y} (u(\mathbf{x}) - \mathcal{N}(\mathbf{x};\theta_s))^2}{\sum_{\mathbf{x} \in Y}  u(\mathbf{x})^2} } \quad \text{and} \quad \epsilon_f (Y) = \sqrt{\frac{\sum_{\mathbf{x} \in Y} ( f(\mathbf{x}) + \triangle \mathcal{N}(\mathbf{x};\theta_s) )^2}{\sum_{\mathbf{x} \in Y}  f(\mathbf{x})^2} },
\end{equation}
respectively. Here, the data set $Y$ is a training data set $X$ used to calculate the training performance or a test data set $\widetilde{X}$ used to calculate the generalization performance.
We note that the training data set $X$ is transformed into the scaled domain $\Omega_b$
when training the neural network parameters $\widehat{\theta}$ for the scaled problem
with the loss function~\eqref{scale:loss:1}.

\subsection{The residual correction method}

In the scaled problem~\eqref{scale:model}, the low-frequency parts of the solution to the original problem become vanishingly small, and the neural network solution $\mathcal{N}(\widehat{\mathbf{x}};\widehat{\theta})$ only works for approximating high-frequency parts. Such a low-frequency part solution thus remains as the dominant error in the scaled-back neural network solution $\mathcal{N}(\mathbf{x};\theta_s)$.
To capture such a low-frequency part of the solution, we include a post-processing step
in which we solve the following residual error equation in the original
problem domain $\Omega$,
\begin{equation}\label{pb:residual-correction}
\begin{split}
-\triangle w(\mathbf{x}) &= f(\mathbf{x}) + \triangle \mathcal{N}(\mathbf{x};\theta_s), \quad \mathbf{x} \in \Omega,\\
w(\mathbf{x})&=g(\mathbf{x})-\mathcal{N}(\mathbf{x};\theta_s), \quad \mathbf{x} \in \partial \Omega.
\end{split}
\end{equation}
In the residual error equation, $\mathcal{N}(\mathbf{x};\theta_s)$ is the scaled-back neural network approximate solution that is obtained using the previous domain scaling method.
The conventional neural network solution $\mathcal{N}(\mathbf{x};\theta)$ can also be used instead of $\mathcal{N}(\mathbf{x};\theta_s)$. 
For the above residual error equation, we find a neural network solution $\mathcal{N}_r(x;\theta_r)$, whose parameters $\theta_r$ are trained for the corresponding loss function,
\begin{equation}\label{residual:loss}
\begin{split}
L_{m,r}(\theta_r) &=w_1 \left( \frac{1}{|X_{\Omega}|} \sum_{\mathbf{x} \in X_{\Omega} } \left| f(\mathbf{x})+\triangle \mathcal{N}(\mathbf{x};\theta_s)+\triangle \mathcal{N}_r(\mathbf{x};\theta_r) \right|^2 \right) \\
& \quad + w_2 \left( \frac{1}{|X_{\partial \Omega}|} \sum_{ \mathbf{x} \in X_{\partial \Omega}} \left| g(\mathbf{x})-\mathcal{N}(\mathbf{x};\theta_s)- \mathcal{N}_r(\mathbf{x};\theta_r) \right|^2 \right).
\end{split}
\end{equation}
In the above loss function, we set the weight factors similarly as in the scaled model problem,
$$w_1=\dfrac{1}{\Vert f+\triangle \mathcal{N}\Vert_{X_\Omega}^2} \quad \text{and} \quad w_2=1.$$
We note that the parameter training in the above residual error problem can be performed efficiently with a relatively small number of training epochs, since the dominant part of the solution to the residual error problem has smooth and low-frequency features that can be well approximated at the early training stage with a relatively small network.
After the post-processing step, we can obtain the neural network approximate solution as the sum of
the scaled-back solution and the residual correction solution,
$$U(\mathbf{x};\theta_s,\theta_r)=\mathcal{N}(\mathbf{x};\theta_s)+\mathcal{N}_r(\mathbf{x};\theta_r).$$

We now introduce some notation that we will use to estimate the training and generalization performance of the residual correction step. We define the relative error $\epsilon_u^r(Y)$ and the relative residual error $\epsilon_f^r(Y)$ in the $l^2$-norm as
\begin{equation}\label{error:residual}
\begin{split}
\epsilon_u^r (Y)= \sqrt{\frac{\sum_{\mathbf{x} \in Y} \left( u(\mathbf{x}) - U(\mathbf{x};\theta_s,\theta_r) \right) ^2}{\sum_{\mathbf{x} \in Y}  u(\mathbf{x})^2} }
\quad \text{and} \quad
\epsilon_f^r (Y)= \sqrt{\frac{\sum_{\mathbf{x} \in Y} \left( f(\mathbf{x}) + \triangle U(\mathbf{x};\theta_s,\theta_r) \right)^2}{\sum_{\mathbf{x} \in Y}  f(\mathbf{x})^2} },
\end{split}
\end{equation}
where the data set $Y$ is a training data set $X$ or a test data set $\widetilde{X}$ to
calculate the training performance and the generalization performance, respectively.

\section{Numerical experiments}\label{sec4}
We present the results of some numerical experiments to confirm the efficiency and accuracy of the proposed method
for solving the multi-frequency model problems.
First, we consider a regression problem for a model with two different frequencies as a preliminary test to allow us to choose appropriate activation functions for this study and ensure that the domain scaling method effectively captures high-frequency components.
We then consider Poisson problems with multi-frequency solutions.

We use the Adam optimizer~\cite{kingma2014adam} to train the network parameters
to the corresponding loss function and set the learning rate as $0.001$.
In one-dimensional cases, we employ a fully connected neural network of $L=4$ hidden layers, with $N=20$ nodes per hidden layer, for both the scaled problem and the residual error problem.
In two-dimensional cases, for the scaled problem we employ a fully connected neural network of $L=4$ hidden layers, with $N=70$ or $100$, and for the residual error problem we employ
a smaller fully connected neural network of $L=4$ with $N=50$.
The network parameters are initialized with the Xavier normal distribution with all initial biases set to 0.01.

For the weight factors $w_1$ and $w_2$ of the loss functions,
we set $w_1$ similarly as in \eqref{weight:w1} and set $w_2$ to one.
The test data set $\widetilde{X}$ is taken to be a set of uniformly spaced points in the problem domain $\Omega$.
We choose 100,000 uniform points in the one-dimensional models and 1,000,000 uniform points in the two-dimensional models as the test data set $\widetilde{X}$.
The training data set $X$ will be provided later depending on the model problem.
For each model problem, we choose six different training data sets $X$ and report the error results
as the mean value across those training sets. We use the same training data set
for both the scaled problem and the residual error problem in our computations.

\subsection{Regression problem with high- and low-frequency components}
In this subsection, we consider a simple regression problem to test the training performance of the neural network, its dependence on the choice of activation function, and the ability of the scaling approach to capture the high-frequency component of the target function.

\begin{figure}[htb!]
 \centering
  \includegraphics[width=8cm]{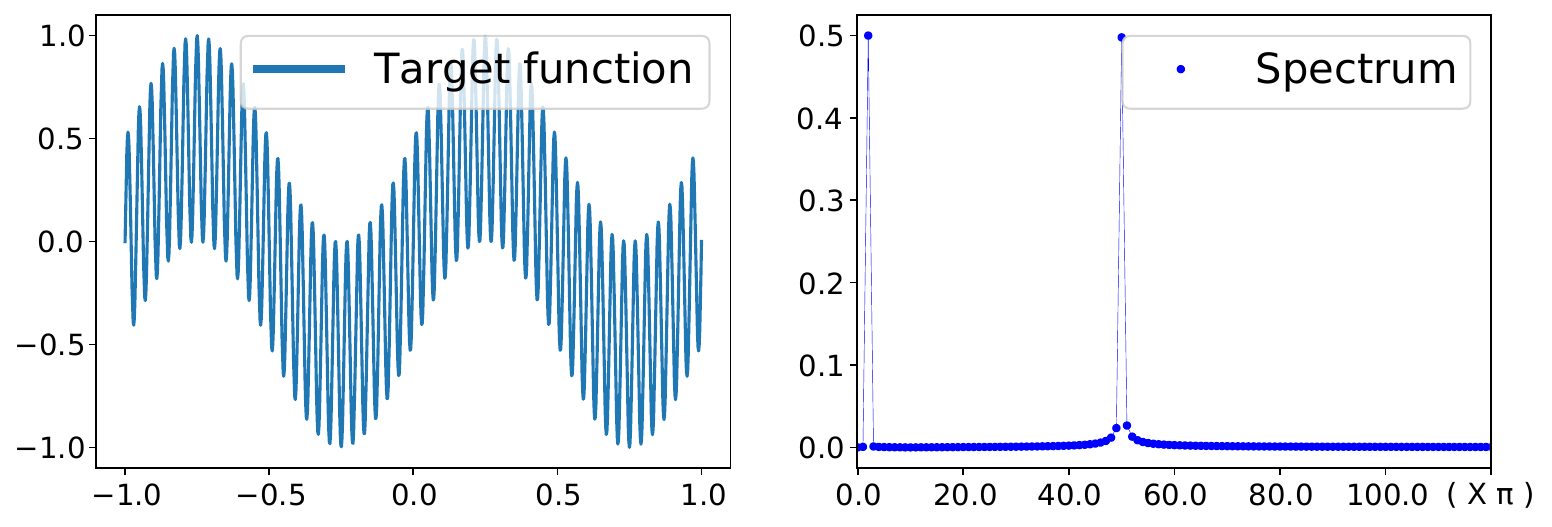}
  \caption{The target function (left) and its spectrum (right) for the regression problem. in \eqref{eq_reg}.}
  \label{sec1_fig1}
  \end{figure}

The target function in our experiment is given as follows:
\begin{equation}\label{eq_reg}
u(x) = \frac{1}{2} \sin (2 \pi x) + \frac{1}{2} \sin (50 \pi x), \quad x \in [-1,\,1].
\end{equation}
Its shape and spectrum are plotted in Fig.~\ref{sec1_fig1}.
  
To find a neural network approximation $\mathcal{N}(x;\theta)$ to the target function $u(x)$, we set the loss function as
\begin{equation*}
L_u (\theta) = \frac{1}{|X|} \sum_{x_i \in X} \left( \mathcal{N}(x_i; \theta) -u(x_i) \right)^2,
\end{equation*}
where $X$ denotes the set of training data points chosen from the problem domain.

We consider six different training sets for $X$: one with 1,000 uniform points and the others with 1,000 randomly selected non-uniform (Latin hypercube) points. We report the mean and standard deviation values of $\epsilon_u(X)$ and $\epsilon_u(\widetilde{X})$ for the six training sets.
We note that the set $\widetilde{X}$ denotes the test data set that is used
to measure the generalization performance of the neural network solution.

We employ various activation functions: 
$\tanh(x)$, Mish, and $\sin(ax)$.
The $\tanh(x)$ activation function is commonly used in regression problems, the Mish function was introduced to overcome restrictions in the ReLU activation function, i.e., it is a differentiable function that approximates the ReLU function~\cite{misra2019mish}, and the $\sin(x)$ activation function is known to be more suitable for approximating derivative information and periodic models~\cite{sitzmann2020implicit}.
We note that a $\sin(ax)$ activation function with a large value $a$ can handle the high-frequency part of the target function more effectively, as can be seen in the numerical results in Table~\ref{sec1_tb1}, and
we thus included this in our numerical computation. 
For the scaled regression model problem, we only employ the $\sin(x)$ activation function and test the scaling method with two different scaling factors, $b=2\pi$ and $50\pi$.

\begin{figure}[htb!]
 \centering
  \includegraphics[width=5cm]{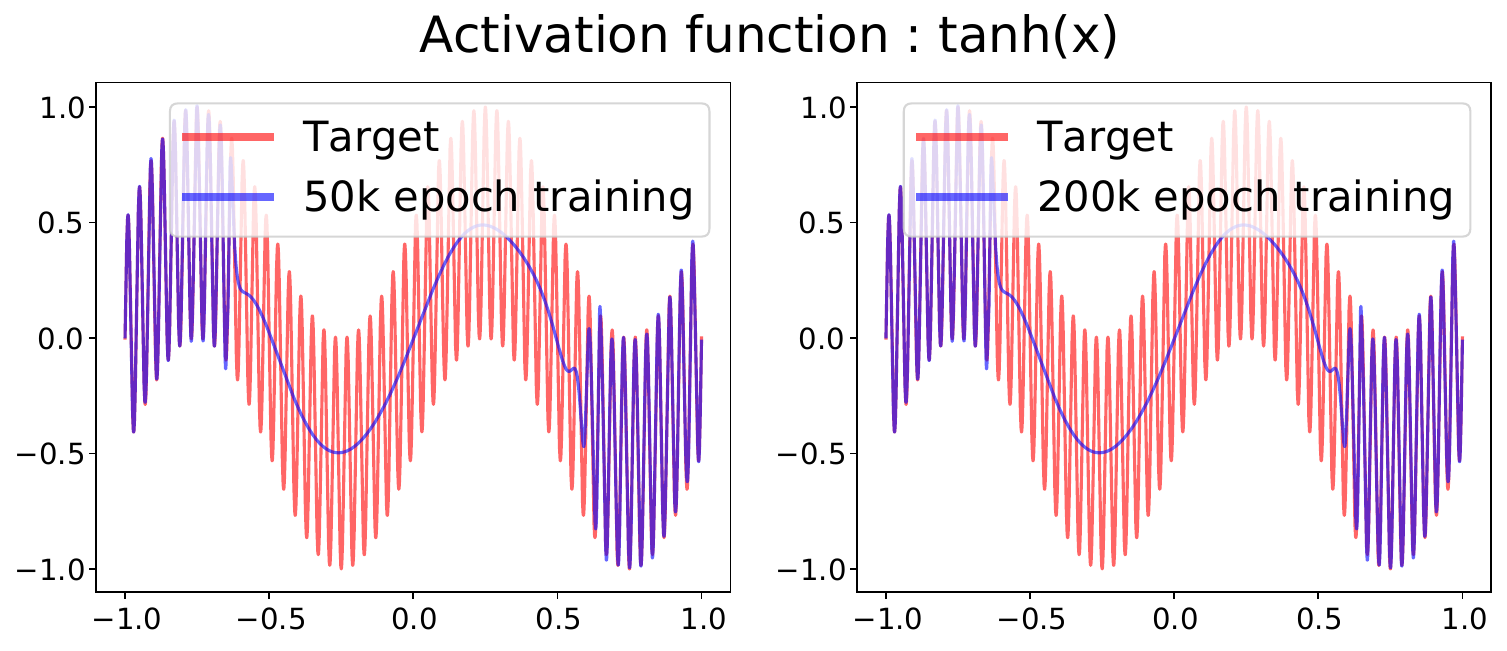} \ \vline \
  \includegraphics[width=5cm]{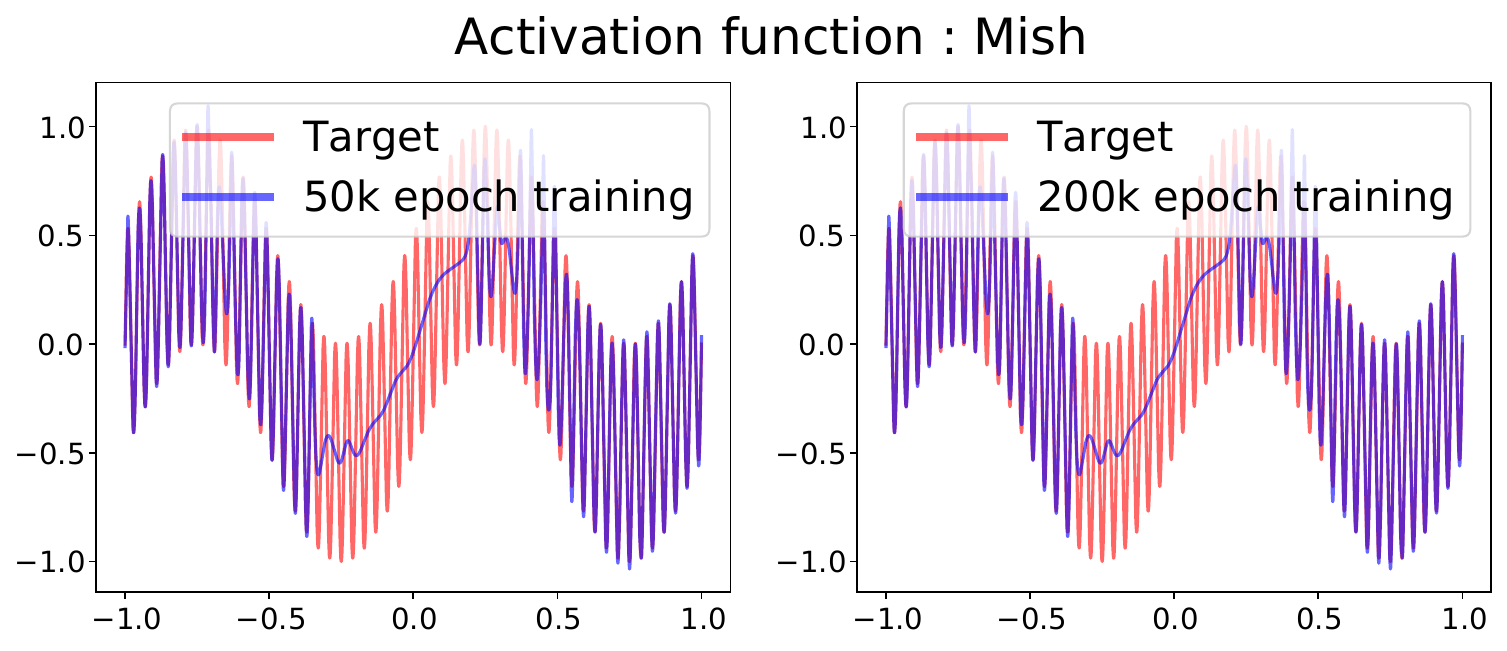} \ \vline \
  \includegraphics[width=5cm]{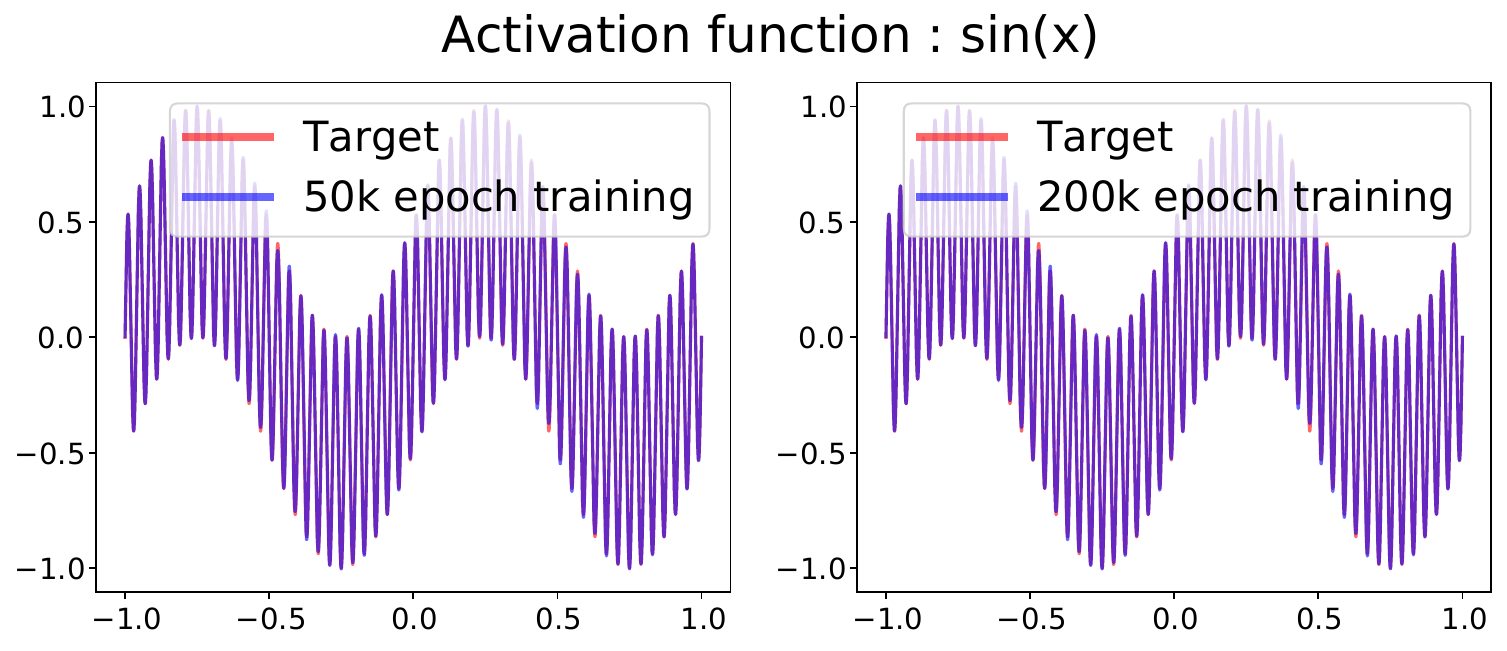}    \\
  \includegraphics[width=5cm]{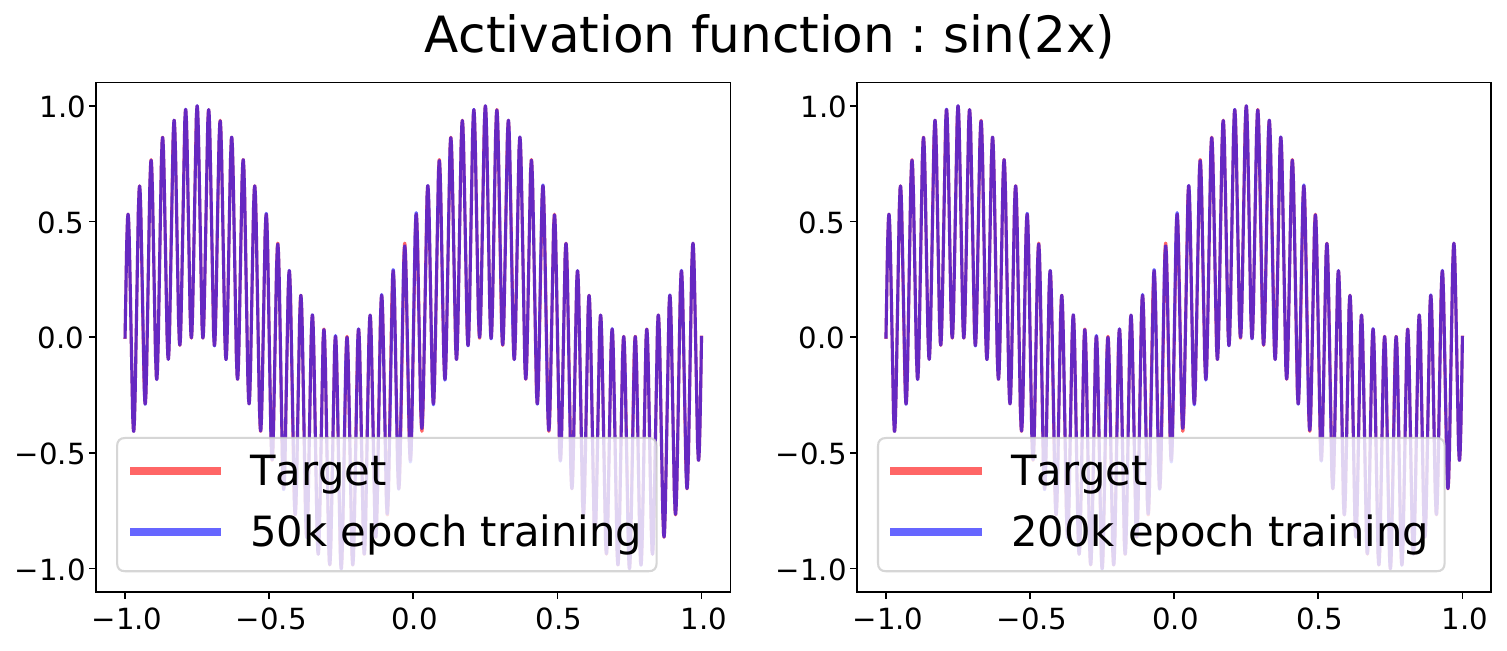}  \ \vline \
  \includegraphics[width=5cm]{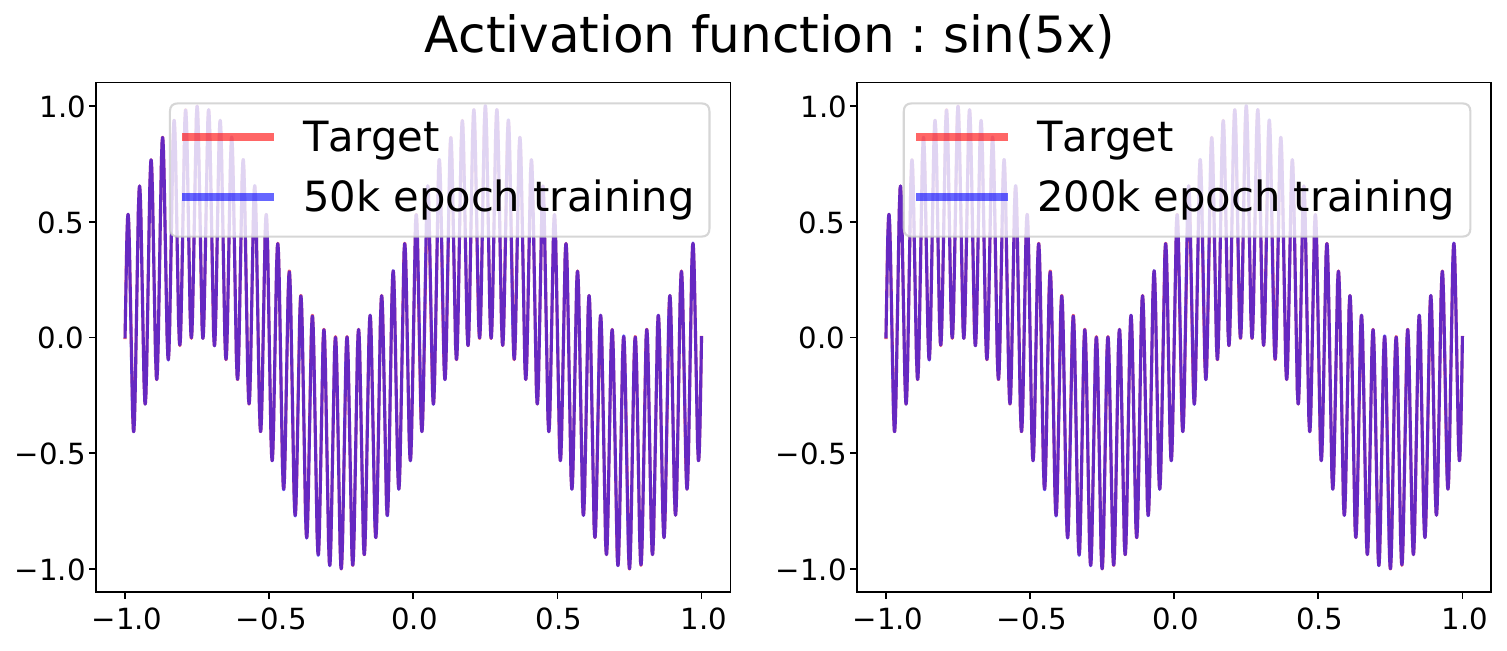}  \ \vline \
  \includegraphics[width=5cm]{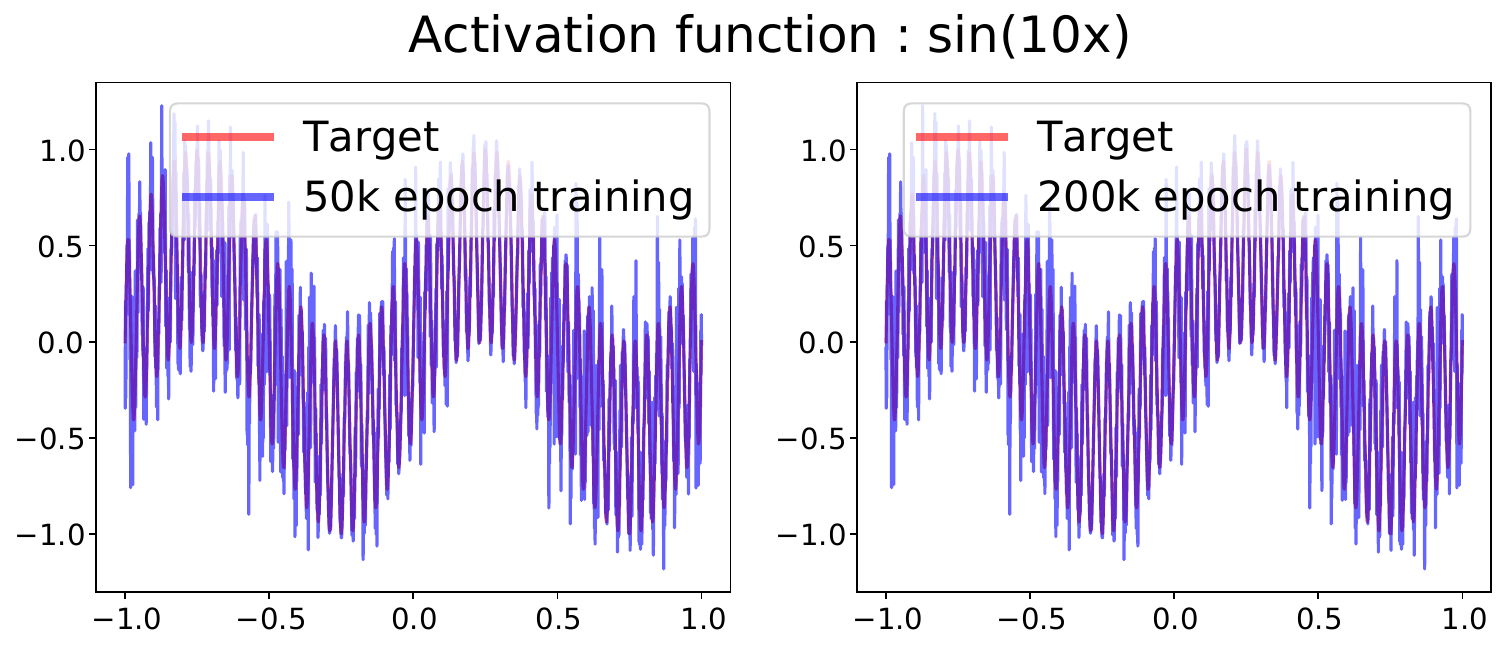} \\ \hspace{-5.35cm}
  \includegraphics[width=5cm]{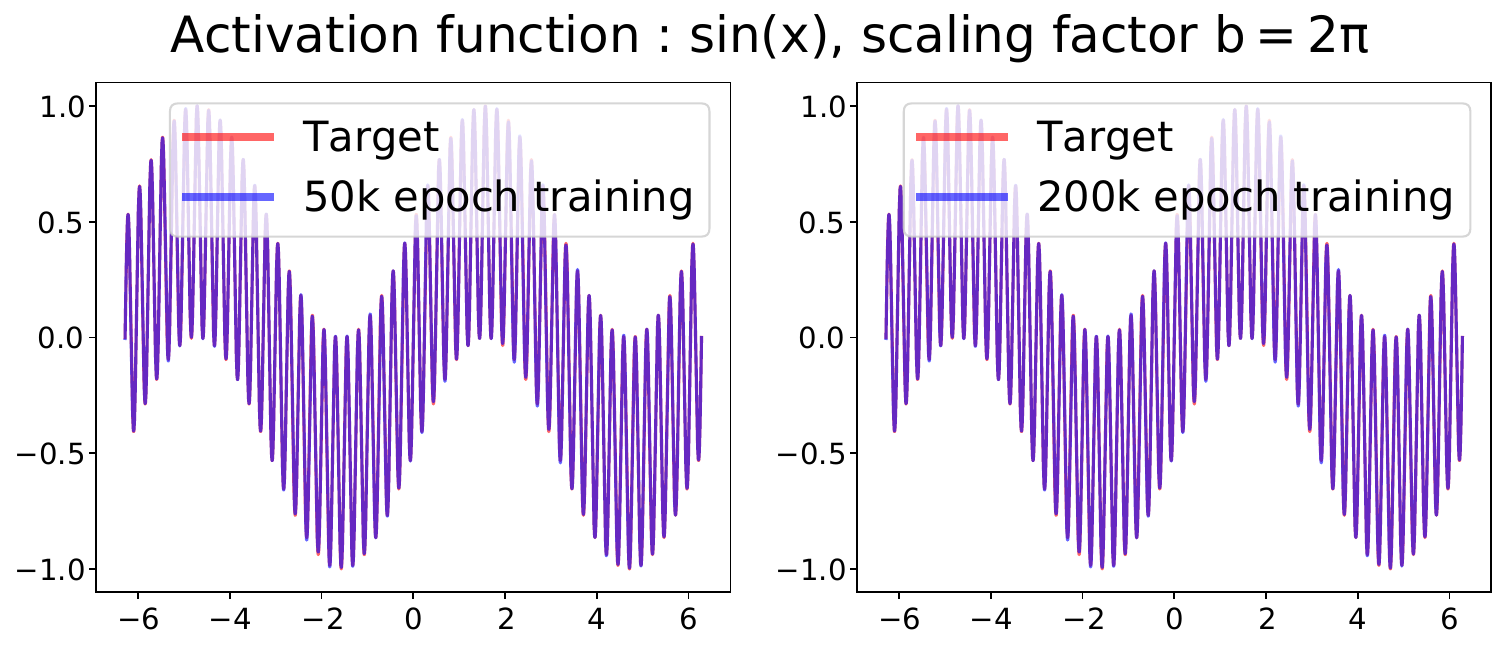} \ \vline \
  \includegraphics[width=5cm]{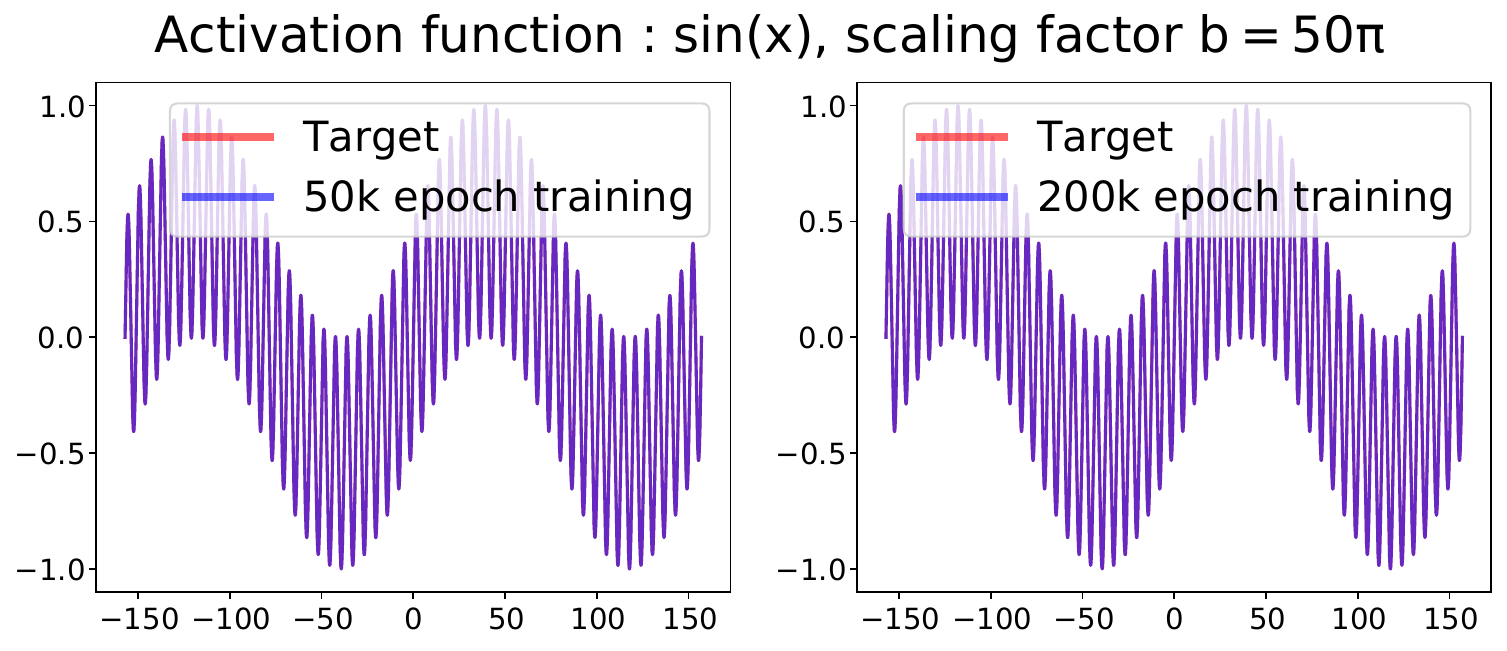} \\
  \caption{Solution plots after 50,000 and 200,000 training epochs for the regression model in \eqref{eq_reg} with various activation functions,
  and for the corresponding scaled regression model with scaling factors $b=2\pi,\; 50\pi$ and with the
  $\sin(x)$ activation function.}
  \label{sec1_fig2}
  \end{figure}

The results from the regression model with various activation functions, as well as the scaled regression model, are presented in Table~\ref{sec1_tb1} and Fig.~\ref{sec1_fig2}.
The $\sin(5x)$ activation function shows the best training and generalization performance. Even at the early training stage of 50,000 epochs, the high-frequency part is well resolved, and the $\sin(5x)$ activation function is thus well suited to approximating the high-frequency solution part.
A choice of $\sin(10x)$ does not result in accurate regression of the target function. As we can see in  Fig.~\ref{sec1_fig2}, there appears to be overfitting and a strong oscillation between the training samples.
The $\tanh(x)$ and Mish activation functions fails to approximate the high-frequency solution part even after 200,000 training epochs.
For the $\sin(x)$ and $\sin(2x)$, the high-frequency solution part is well resolved after 200,000 training epochs.
For the regression problem with both high- and low-frequency parts, the choice of activation function greatly affects the training performance and accuracy.
The $\sin(ax)$ activation functions are better than the $\tanh(x)$ and Mish activation functions. For the scaled regression model with $b=50\pi$, we observe that it shows
the best training and generalization performance.
The errors are even smaller than those with the $\sin(5x)$ activation function.

\begin{table}[ht!]
\begin{center}
  \caption{Error values $\epsilon_u(X)$ and $\epsilon_u(\widetilde{X})$ after 200,000 training epochs for the regression model in \eqref{eq_reg} with various activation functions
  and for the corresponding scaled regression model with scaling factors $b=2\pi,\; 50\pi$ and with the
  $\sin(x)$ activation function
  The numbers listed are the means and standard deviations of six different training sample sets $X$.
  The results are extracted from the spot with the smallest loss value during the training epochs.
 }\label{sec1_tb1}
{\footnotesize \renewcommand{\arraystretch}{1.2}
    \begin{tabular}{ccc|ccc}
        \Xhline{3\arrayrulewidth}
  Activation function & ${\epsilon}_u(X)$ & $\epsilon_u(\widetilde{X})$ & Activation function & ${\epsilon}_u(X)$ & $\epsilon_u(\widetilde{X})$
\\[0.3mm]         \Xhline{0.8\arrayrulewidth}
$\tanh(x)$             & 5.37e-01 {\scriptsize(6.69e-02)}& 5.38e-01 {\scriptsize(6.64e-02)}&
$\sin(5x)$               & 2.07e-03 {\scriptsize(4.20e-04)}&2.11e-03 {\scriptsize(4.28e-04)} \\
Mish                       & 4.91e-01 {\scriptsize(5.64e-02)}& 4.91e-01 {\scriptsize(5.63e-02)}&
$\sin(10x)$            & 2.53e-02 {\scriptsize(2.52e-02)}& 1.82e-01 {\scriptsize(2.01e-01)}  \\
$\sin(x)$      & 8.74e-02 {\scriptsize(1.77e-01)}& 8.76e-02 {\scriptsize(1.77e-01)}&
$\sin(x), \ b=2 \pi$ & 9.90e-03 {\scriptsize(8.38e-03)}& 9.97e-03 {\scriptsize(8.46e-03)} \\
$\sin(2x)$               & 5.97e-03 {\scriptsize(1.28e-03)}& 5.99e-03 {\scriptsize(1.28e-03)}&
$\sin(x), \ b=50\pi$ & 5.06e-04 {\scriptsize(4.01e-05)}&5.38e-04 {\scriptsize(4.32e-05)}\\
[0.5mm]     \Xhline{3\arrayrulewidth}
    \end{tabular}
}

\vskip-.7truecm
\end{center}
\end{table}

In addition, we define $\alpha_u^l$ and $\alpha_u^h$ as the amplitudes of the low- and high-frequency errors between the target and trained neural network, $\{r_j := u(x_j)-\mathcal{N}(x_j,\theta) \ | \ x_j \in X \}$, respectively, where those values are calculated by applying the discrete Fourier transform,  
$F_k = \sum_{j=0}^{N-1} r_j  e^{- \frac{i 2 \pi}{N} kj}$ with $N=|X|$.
For the target model \eqref{eq_reg}, we calculate $\alpha_u^l = F_{2}$ and $\alpha_u^h = F_{50}$.

\begin{figure}[htb!]
 \centering
  \includegraphics[width=3.6cm]{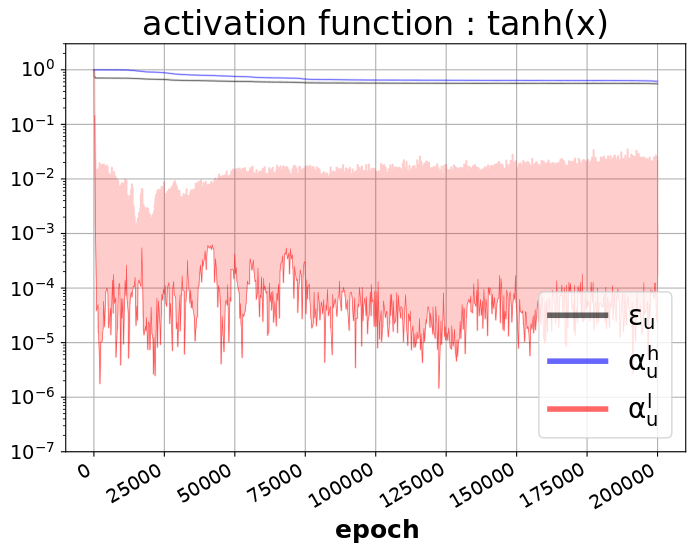}
 \includegraphics[width=3.6cm]{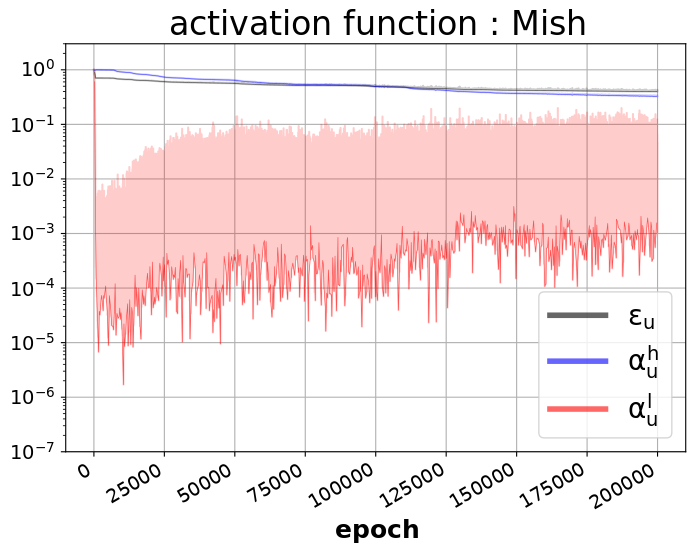}
 \includegraphics[width=3.6cm]{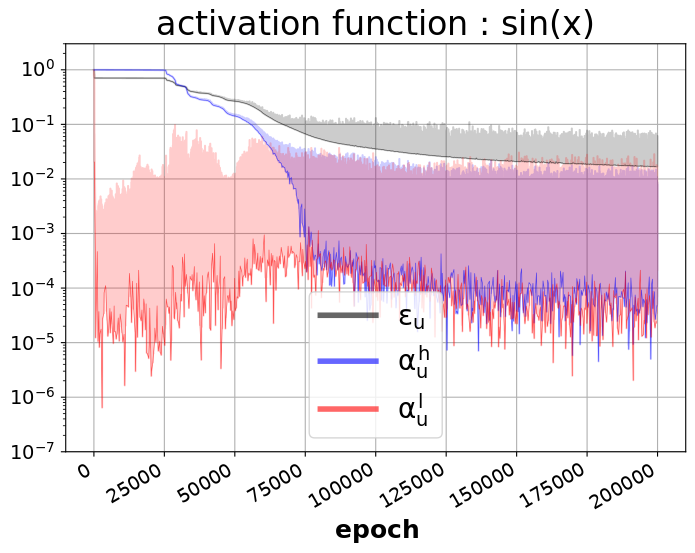}
  \includegraphics[width=3.6cm]{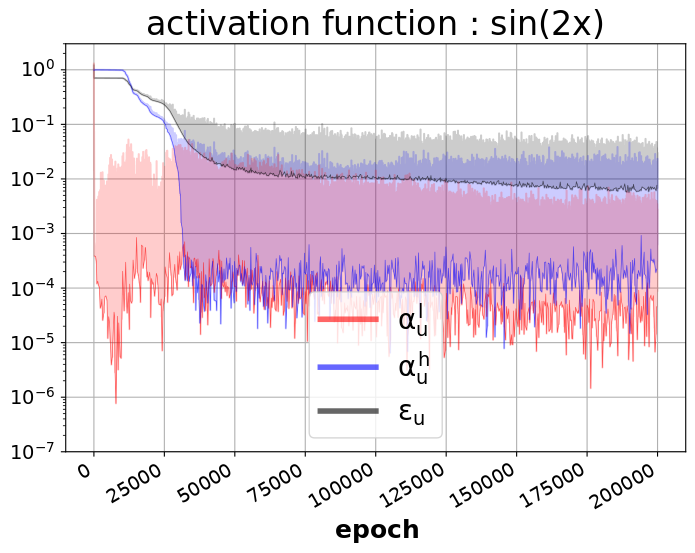}
  \includegraphics[width=3.6cm]{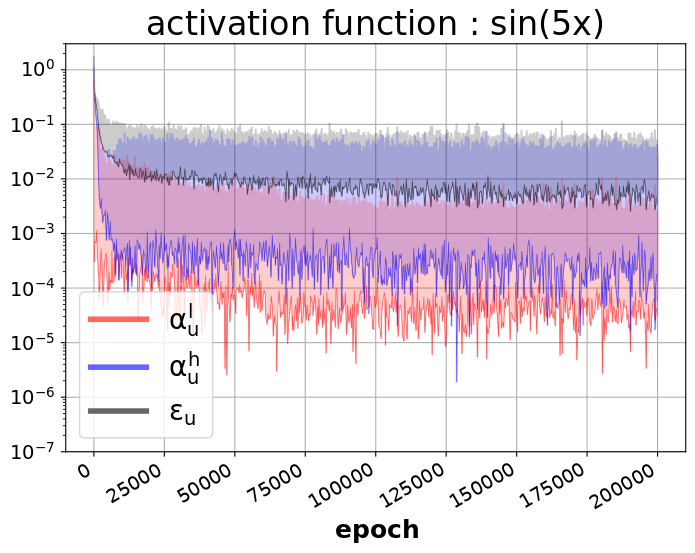}
  \includegraphics[width=3.6cm]{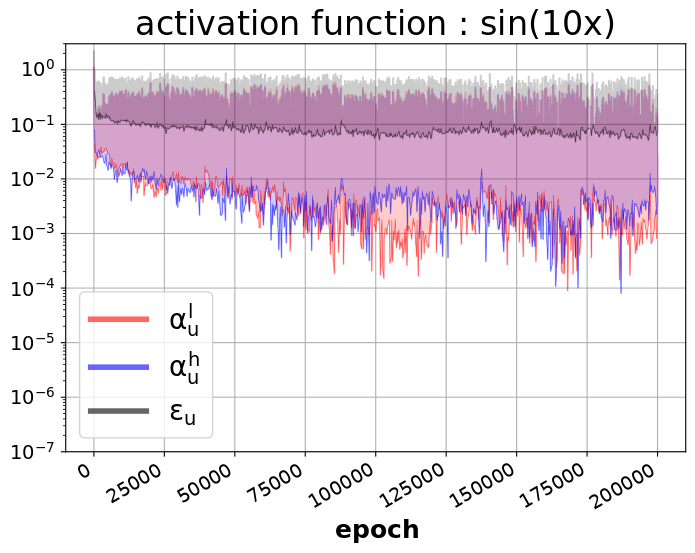}
  \includegraphics[width=3.6cm]{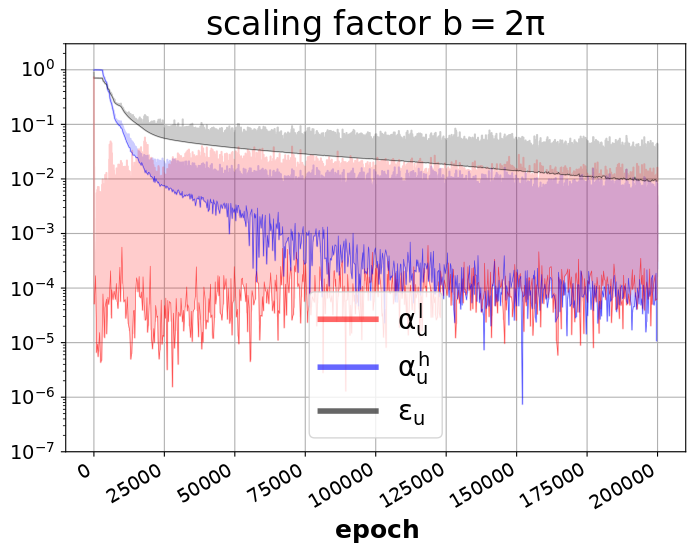}
 \includegraphics[width=3.6cm]{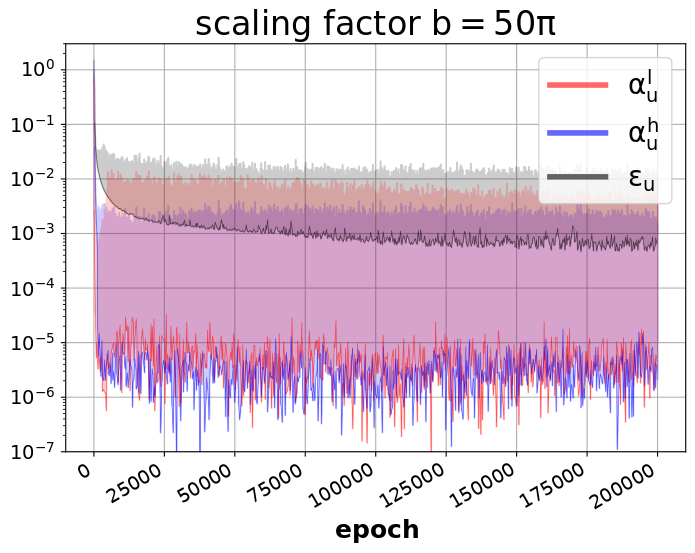}
  \caption{Plots of values $\epsilon_u(X)$, $\alpha_u^h$, and $\alpha_u^l$ over the training epochs for the regression model \eqref{eq_reg} with various activation functions,
  and for the corresponding scaled regression model with scaling factors $b=2\pi,\; 50\pi$ and with the
  $\sin(x)$ activation function.}
  \label{sec1_fig3}
  \end{figure}

In Fig.~\ref{sec1_fig3},
the error values $\epsilon_u$ and the amplitudes of the high- and low-frequency errors,
$\alpha_u^h$, and $\alpha_u^l$, are plotted over the 200,000 training epochs.
The low-frequency amplitude error $\alpha_u^l$, represented by the red lines in the plots,
is sufficiently small during the early training epochs for all activation functions
except the $\sin(10x)$ activation function.
The high-frequency amplitude error $\alpha_u^h$, represented by the blue lines,
is smaller during the early training epochs for the $\sin(ax)$ type activation functions with $a>1$, which means that the higher-frequency activation function is more effective at capturing the high-frequency part of the solution.
As mentioned earlier, the use of excessively high-frequency activation functions may result in overfitting, causing a test error much larger than the training error, as seen in the $\sin(10x)$ activation function.
We also plot the error values $\epsilon_u$, and the amplitudes of high- and low-frequency errors, $\alpha_u^h$, and $\alpha_u^l$, in the scaled regression model with the sine activation function, $\sin(x)$.
In the scaled regression model, we choose $b$ as $2\pi$ and $50\pi$, which are the target function's low- and high-frequency components, respectively.
With the larger scaling factor $b=50 \pi$, the high-frequency part is effectively learned during the early training stage.

From the numerical results with various activation functions, we can observe that the $\sin(ax)$ type activation functions are more suitable for approximating oscillation models with both high- and low-frequency parts. Moreover, we can use the domain scaling approach to enhance the training performance further for approximating the high-frequency component rather than simply using the high-frequency sine activation functions.

\subsection{One-dimensional Poisson problem}
We solve the one-dimensional Poisson problem~\eqref{model} with a multi-frequency solution.
The differential operator amplifies high-frequency parts of the multi-frequency solution, and
the corresponding amplitude for the high-frequency parts in the right-hand side function $f$ can thus become excessively large even when the amplitude of high-frequency parts of the solution $u$ is small.
To demonstrate this phenomenon, we consider the model problem solution as
\begin{equation}\label{poi_test_ft}
u(x) = 5\sin(\pi x) + \sin(8 \pi x) + \frac{1}{2}\sin(16 \pi x) + \frac{1}{4}\sin(32 \pi x) + \frac{1}{8}\sin(64 \pi x).
\end{equation}
in the domain $\Omega = (-1,\; 1)$.

\begin{figure}[htb!]
 \centering
  \includegraphics[width=8cm]{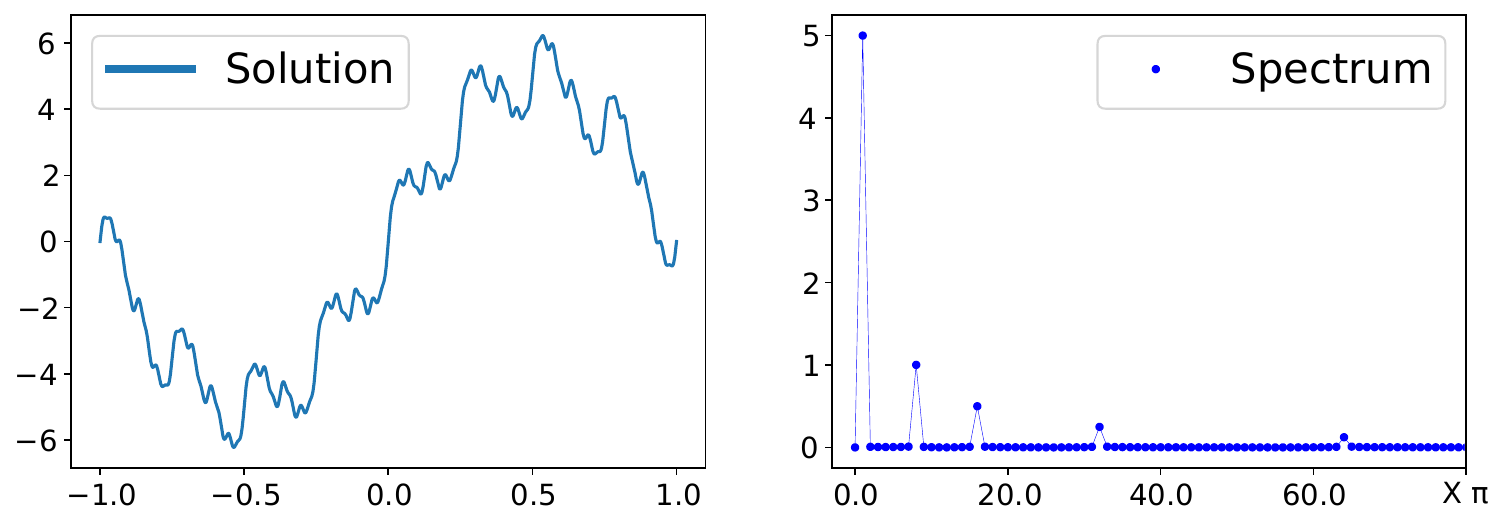}
  \includegraphics[width=8cm]{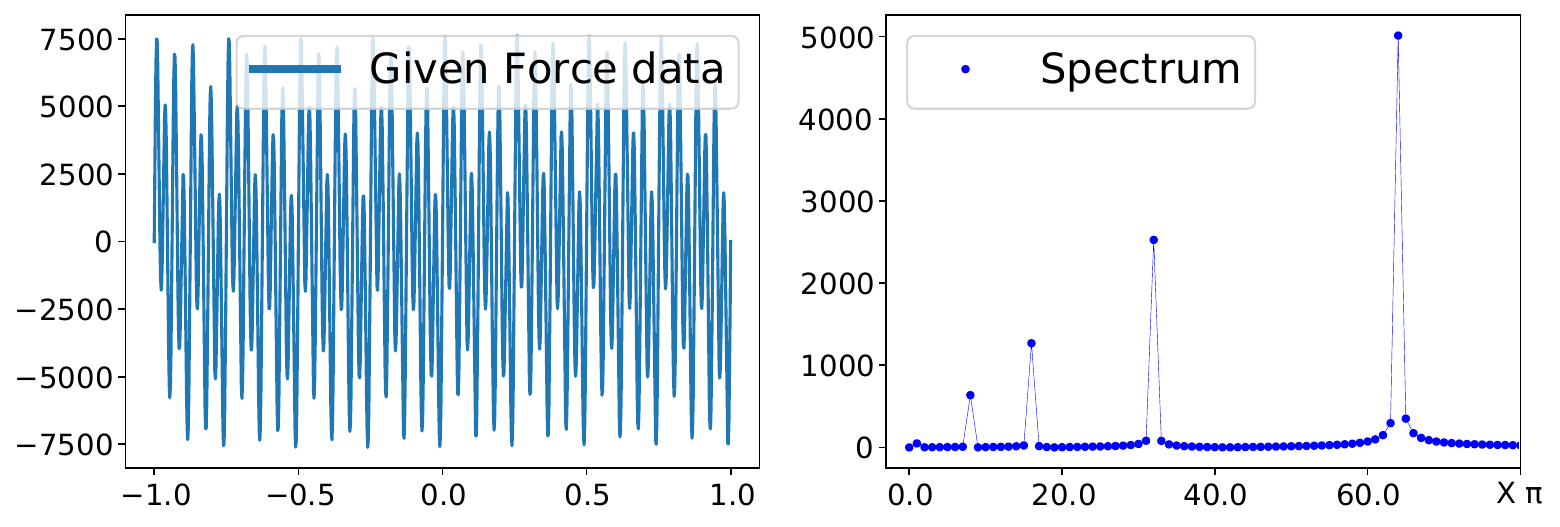}
  \caption{The target solution $u(x)$ and its spectrum (left two graphs) and the right hand side function $f(x)$ and its spectrum (right two graphs).}
  \label{sec2_fig1}
\end{figure}

The shapes and spectra of the solution $u$ and the corresponding right hand side function $f$ in the Poisson model problem are presented in Fig.~\ref{sec2_fig1}.
The small amplitude of the high frequency component in the model solution
appears as the largest amplitude in $f$, as seen in Fig.~\ref{sec2_fig1}.
The Poission model problem with such a magnified high frequency component is known to be difficult to train with
the standard PINN approach \cite{luo2021frequency, xu2019frequency}. We will apply the domain scaling method and residual correction method to the above model problem.

\subsubsection{Activation functions and domain scaling approach}
We examine the neural network training performance for the Poisson problem with various activation functions and for the scaled Poisson problem with various scaling factors. 

\begin{figure}[htb!]
 \centering
  \includegraphics[width=6cm]{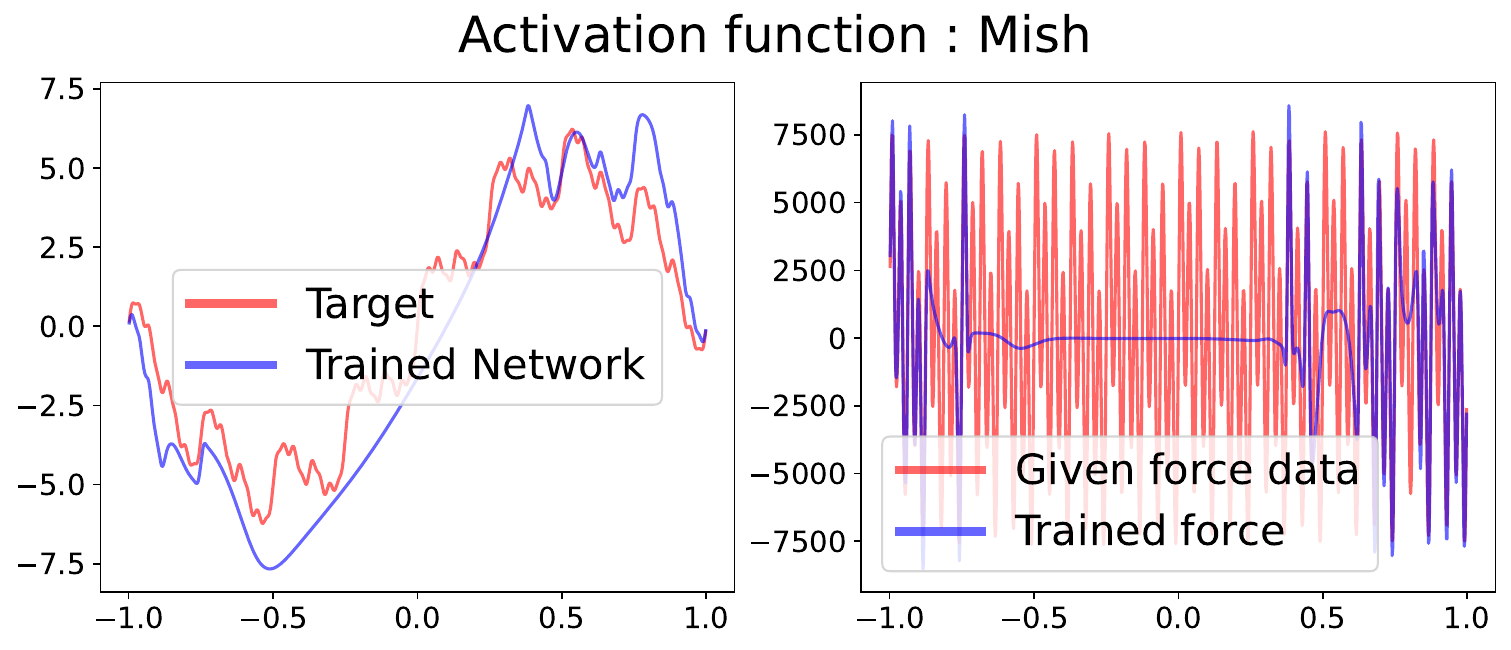}
  \includegraphics[width=6cm]{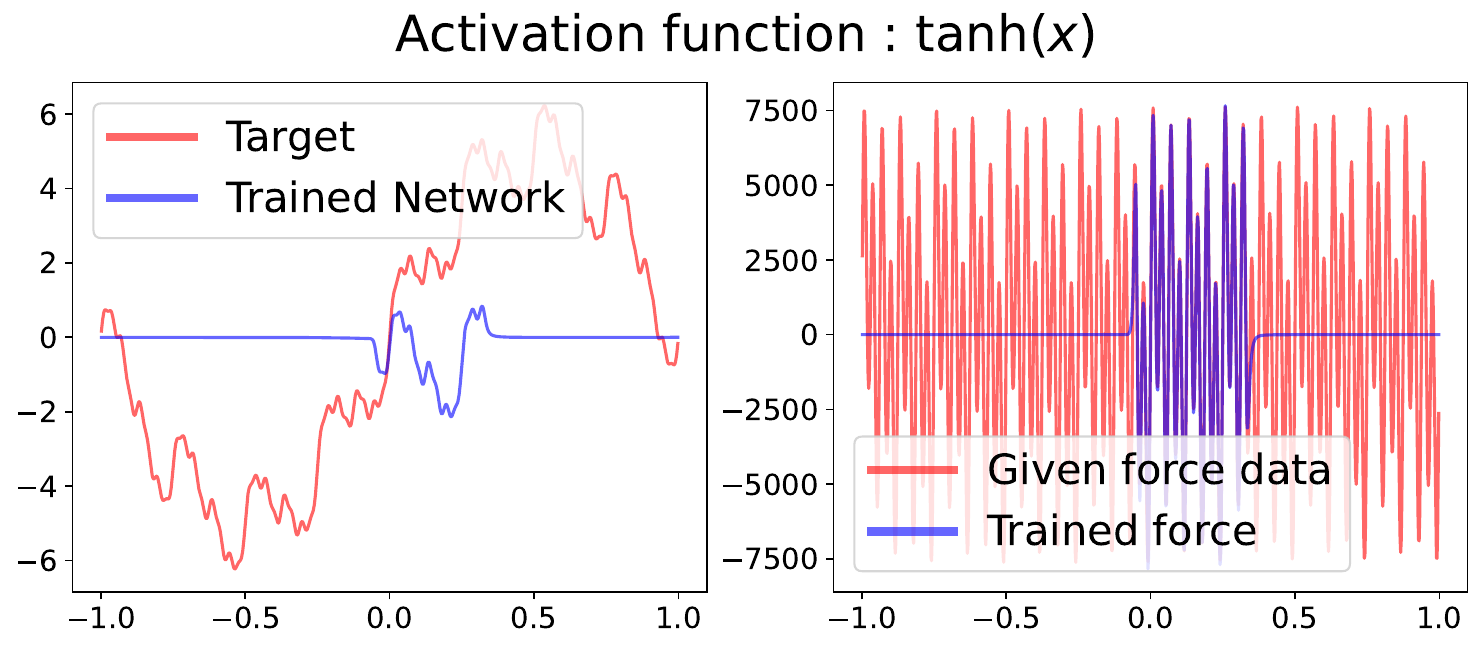} \\
  \includegraphics[width=6cm]{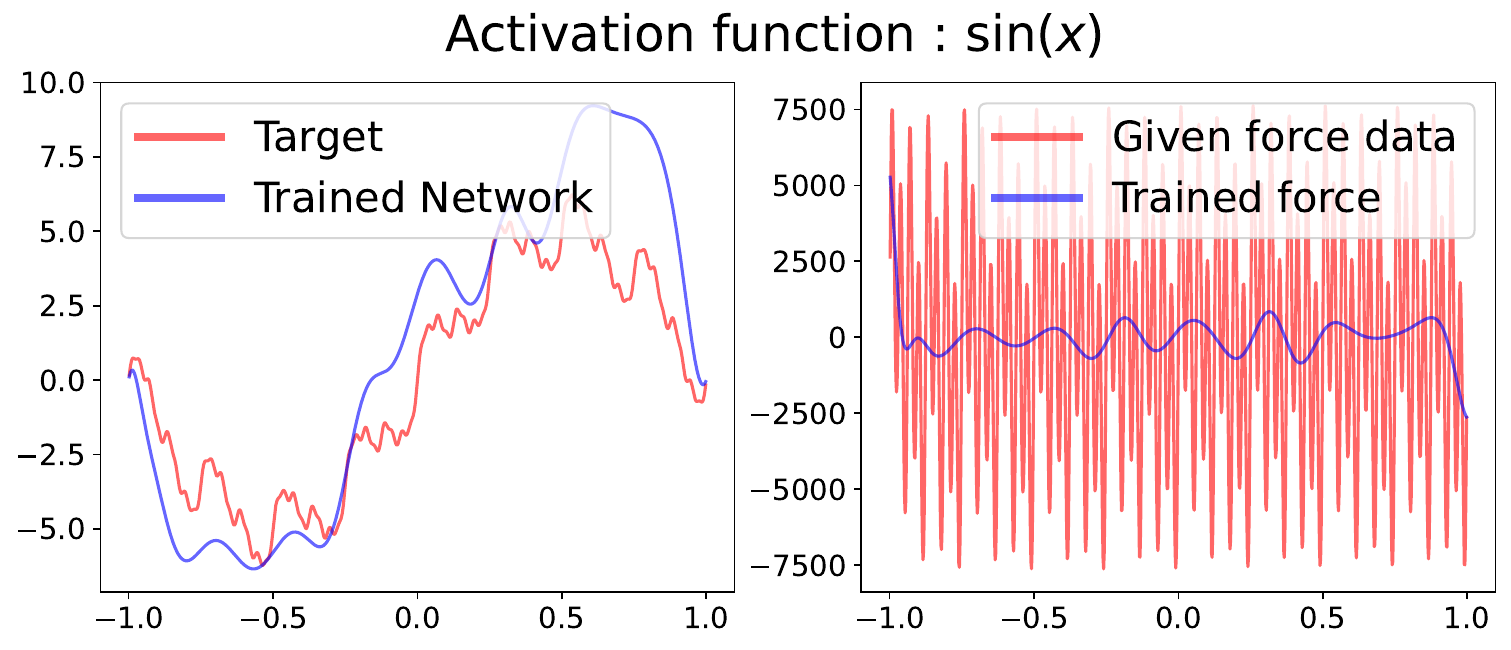}
  \includegraphics[width=6cm]{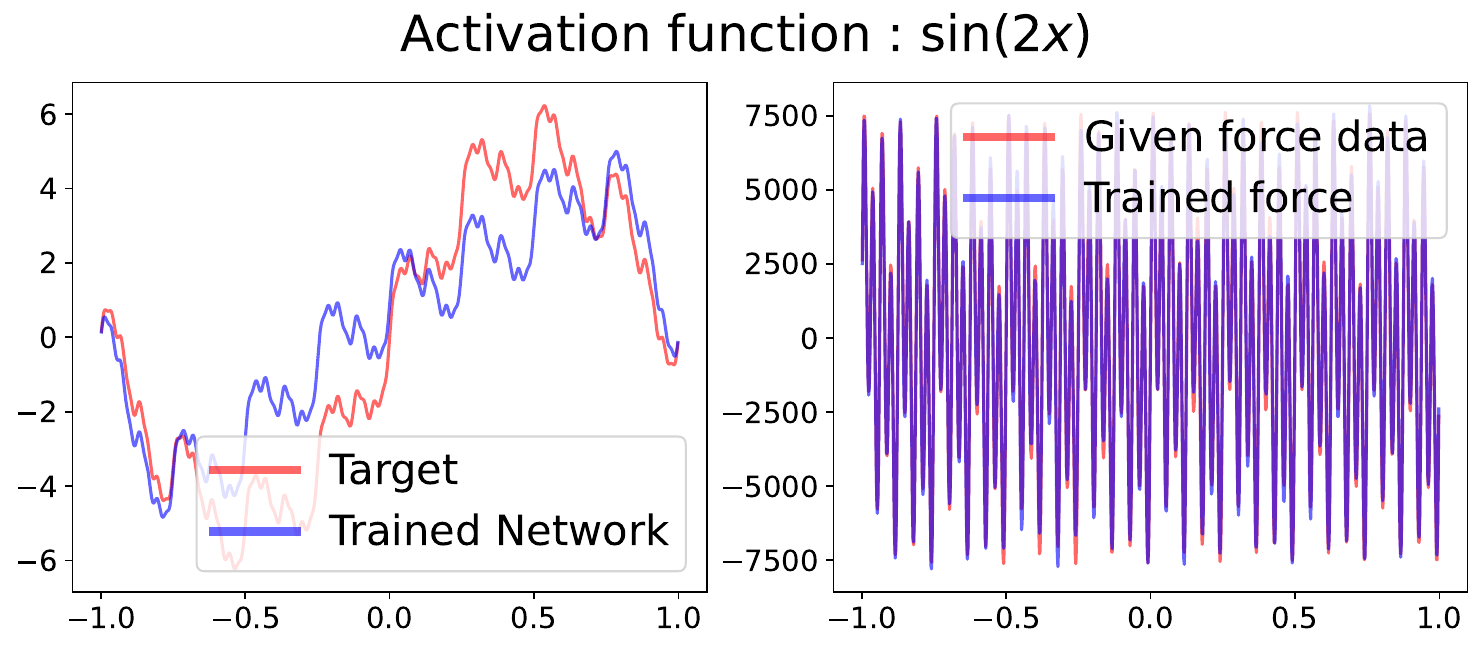} \\
  \includegraphics[width=6cm]{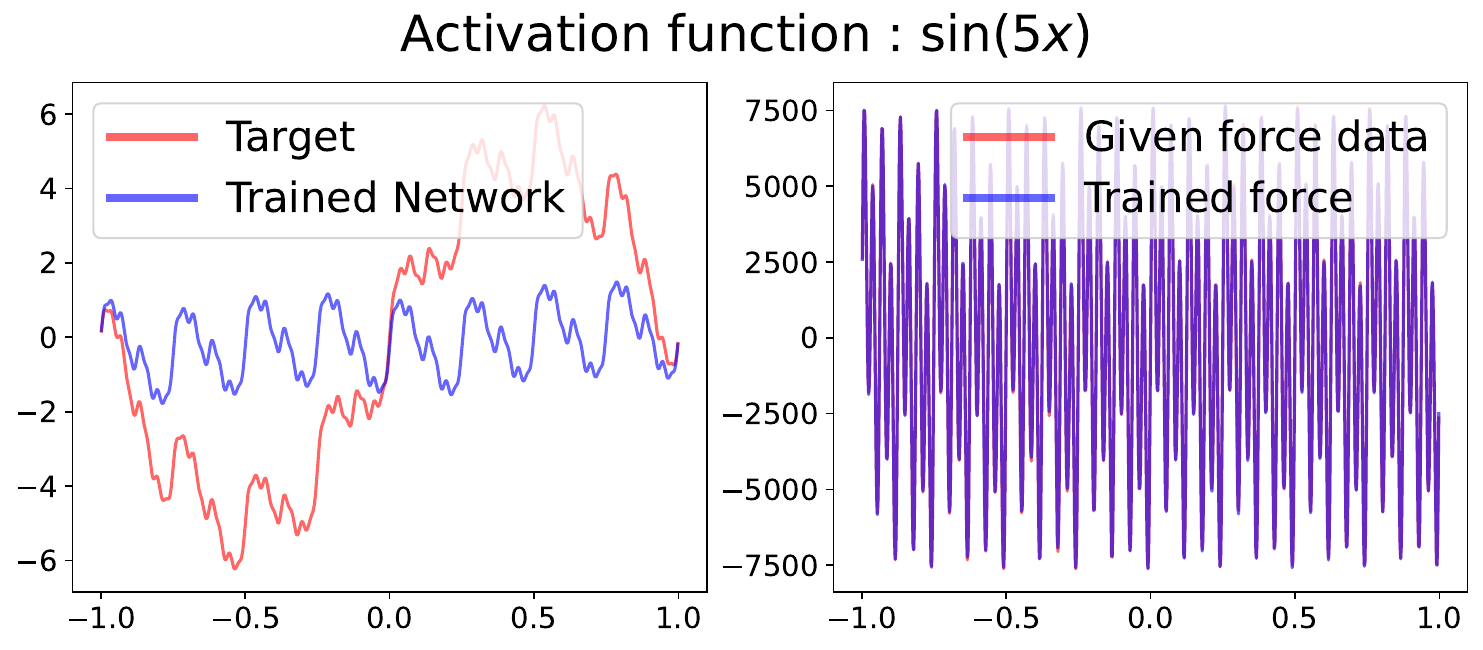}
  \includegraphics[width=6cm]{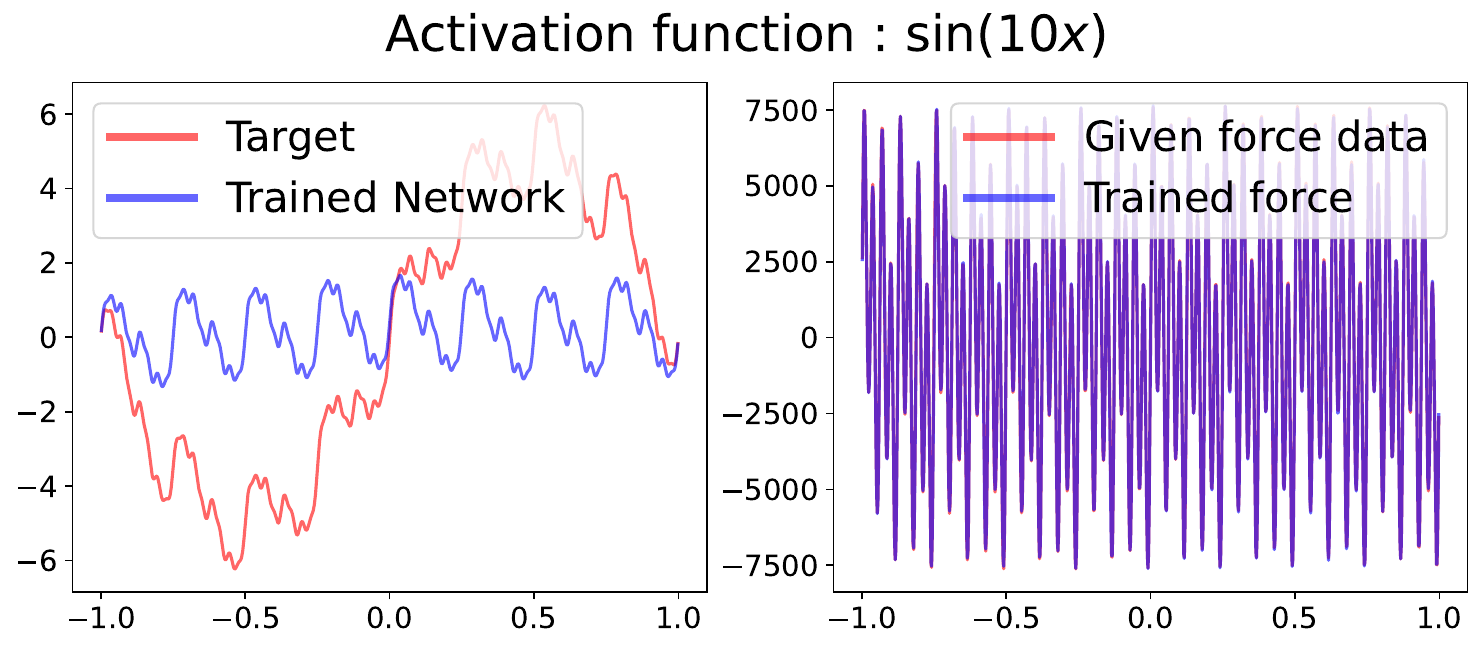} \\
  \includegraphics[width=6cm]{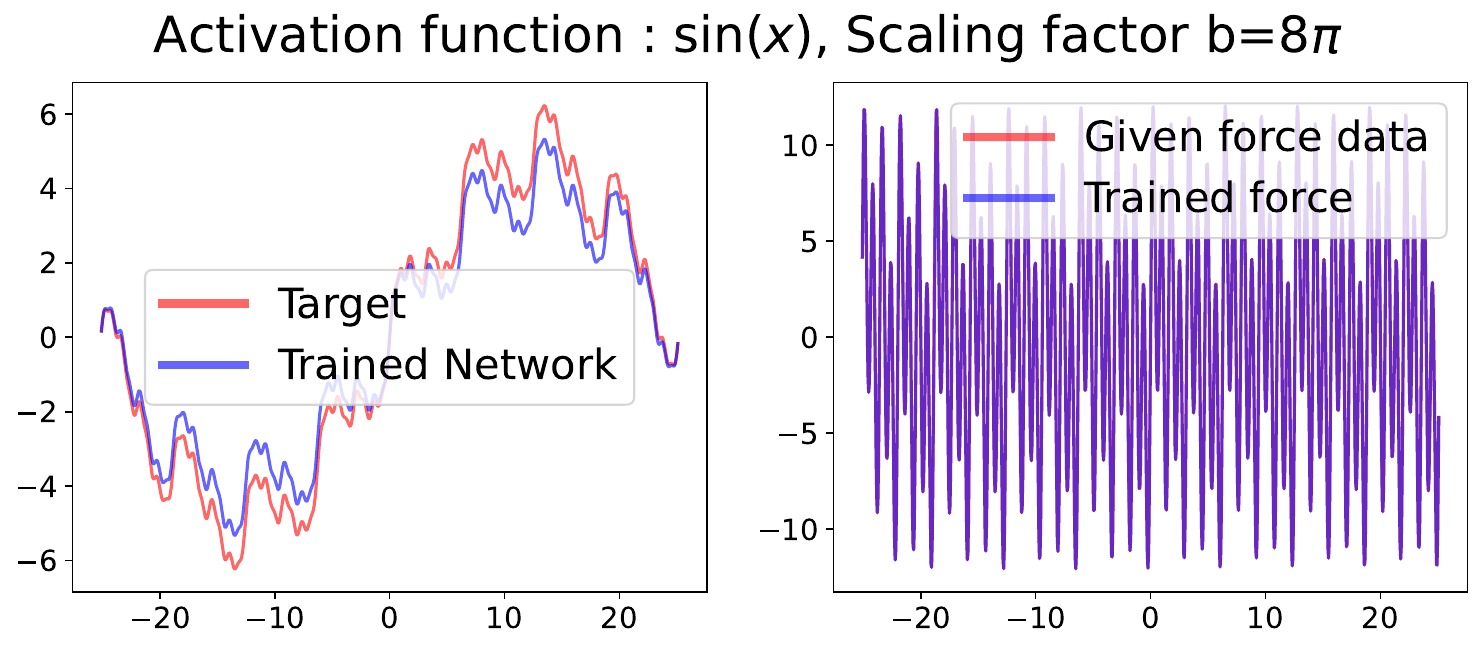}
  \includegraphics[width=6cm]{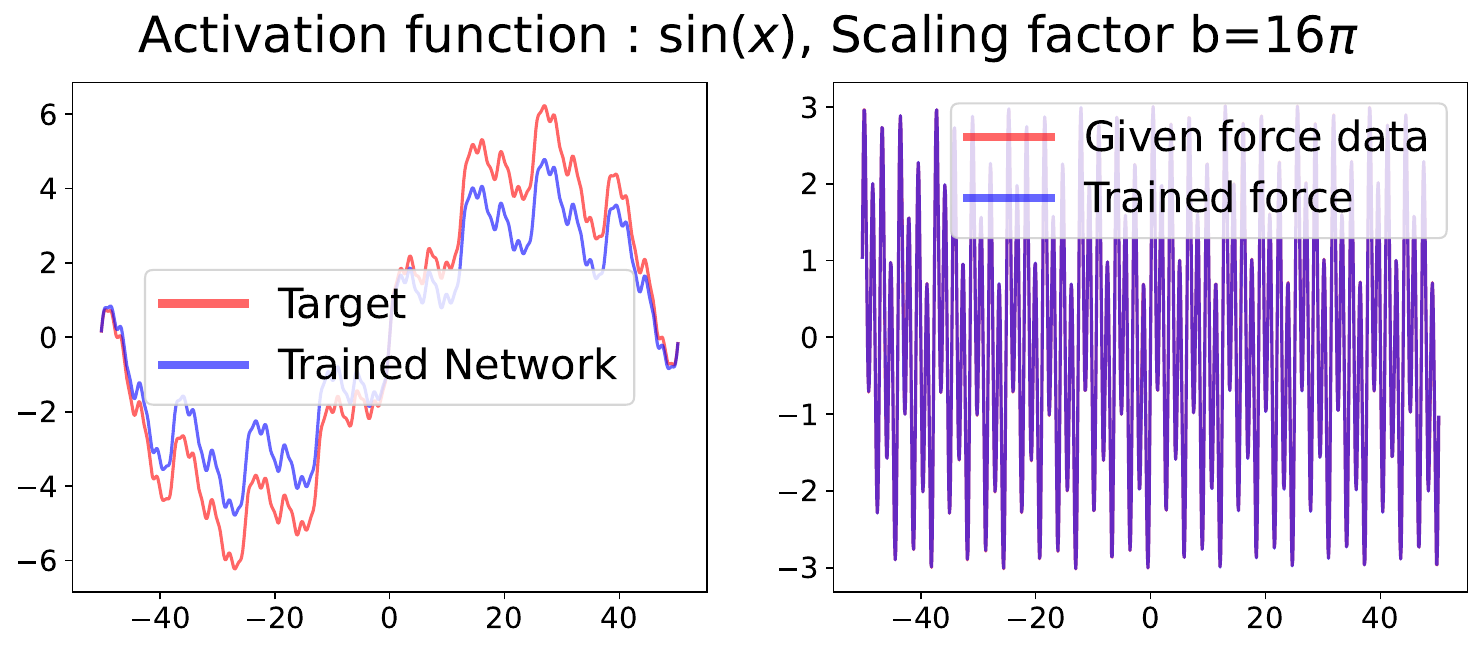} \\
  \includegraphics[width=6cm]{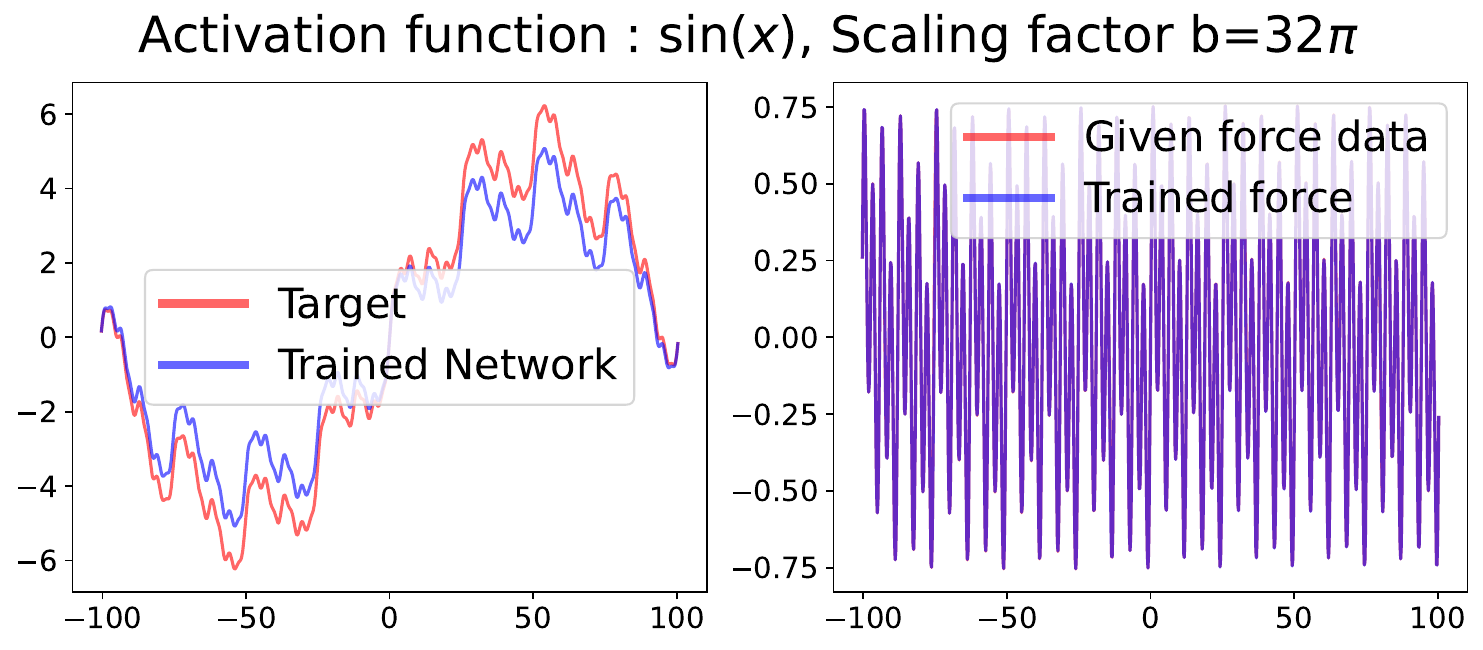}
  \includegraphics[width=6cm]{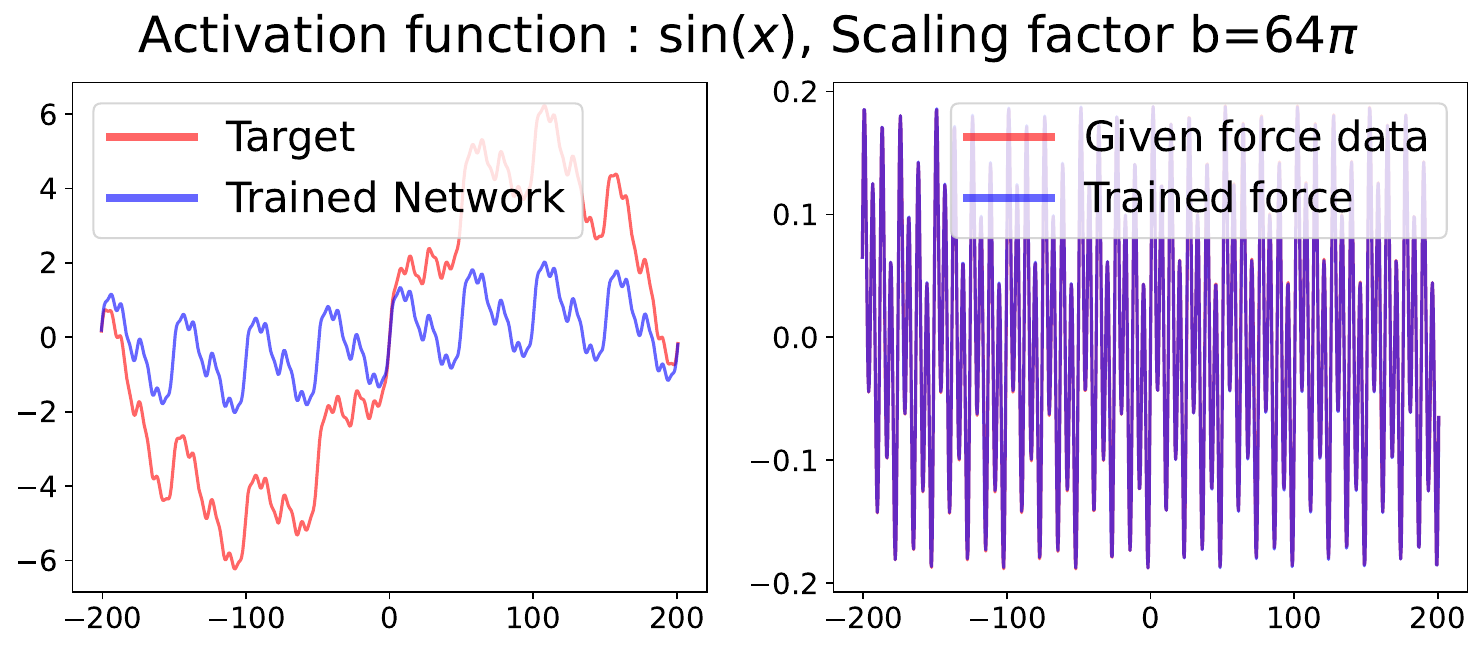} \\
  \caption{Solution and force plots after 100,000 training epochs for the Poisson model~\eqref{model} with the target solution~\eqref{poi_test_ft} using various activation functions
  and for the corresponding scaled Poisson model using various scaling factors with the
  $\sin(x)$ activation function.}
  \label{sec2_fig4}
  \end{figure}

In Fig.~\ref{sec2_fig4}, the target solution $u(x)$ and the trained neural network solution $\mathcal{N}(x;\theta)$, as well as their corresponding right-hand side functions $f(x)$ and $-\triangle \mathcal{N}(x;\theta)$, are plotted for various activation functions: $\tanh(x)$, Mish, and $\sin(ax)$ with $a=1, 2, 5, 10$. 
The results show similar behaviors as in the previous regression model problem.
The $\tanh(x)$ and Mish activation functions could not resolve the high-frequency parts of the solution. With the $\sin(ax)$ activation functions, the high-frequency parts are better approximated.
However, the error of the low-frequency part becomes more dominant as the value of $a$ is increased.

In addition, the target solution $u(x)$ and the trained neural network solution $\mathcal{N}(x;\theta_s)$, as well as their corresponding right-hand side functions $f(x)$ and $-\triangle \mathcal{N}(x;\theta_s)$, are plotted for scaling factors $b=8\pi, 16\pi, 32\pi, 64\pi$. The results for the scaled Poisson problem show that the high-frequency parts of the model solution
are well approximated by the trained neural network $\mathcal{N}(x;\theta_s)$ for all these choices of the scaling factor $b$.

In Table~\ref{sec2_tb1}, we report the quantitative errors for the Poisson problem with the $\sin(ax)$-type activation functions
and for the scaled Poisson problem with various scaling factors $b$.
The $\sin(2x)$ activation function gives the best solution accuracy but the errors are still greater than $0.1 (10\%)$.
For the right-hand side function error $\epsilon_f$, the $\sin(5x)$ activation function gives the best result.
For the scaled Poisson problem, the errors are reduced slightly at the early training epochs, up to 20,000 epochs, and they are further reduced at the later training epochs, up to 100,000 epochs.
We note that in the Poisson problem with the $\sin(5x)$ activation function, we cannot readily observe such an error reduction at the later training epochs as we have seen for the scaled Poisson problem with the scaling factor $b=16 \pi$.

Based on these results, we can conclude that for the multi-frequency Poisson model problem, the simple sine activation function $\sin(x)$  with the domain scaling approach results in better solution accuracy and training efficiency than the $\sin(ax)$-type activation functions for the Poisson problem. 
\begin{table}[h!]
\begin{center}
  \caption{Error values $\epsilon_u(\widetilde{X})$ and $\epsilon_f({X})$ for the Poisson model~\eqref{model} with the target solution~\eqref{poi_test_ft} using various activation functions
  and for the corresponding scaled Poisson model using various scaling factors with the
  $\sin(x)$ activation function.
  The mean values and standard deviations are listed for six different training sample sets $X$.
  The results are extracted from the spot with the smallest loss value during the training epochs.
  }\label{sec2_tb1}
{\footnotesize \renewcommand{\arraystretch}{1.2}
    \begin{tabular}{cccccc}
        \Xhline{3\arrayrulewidth}
& $\sin(x)$ & $\sin(2x)$ & $\sin(5x)$ & $\sin(10x)$ \\  [-0.8ex]
Max.\ epoch & $\epsilon_u(\widetilde{X}), \qquad \ \ \epsilon_f(X)$ & $\epsilon_u(\widetilde{X}), \qquad \ \  \epsilon_f(X)$ &
 $\epsilon_u(\widetilde{X}), \qquad \ \  \epsilon_f(X)$ & $\epsilon_u(\widetilde{X}), \qquad \ \ \epsilon_f(X)$ \\[1mm]         \Xhline{0.8\arrayrulewidth}
20,000
& 1.03e+00 \ \ \ 	1.00e+00
& 9.72e-01 \ \ \ 	6.58e-01
& 9.28e-01 \ \ \ 	2.55e-02
& 1.00e+00 \ \ \ 	5.18e-01  \\ [-1.2ex]
&  \ {\scriptsize(3.64e-02)} \ \ \ {\scriptsize(1.22e-05)}
&  \ {\scriptsize(1.16e-01)} \ \ \ {\scriptsize(1.36e-01)}
& \ {\scriptsize(3.03e-02)} \ \ \ {\scriptsize(5.71e-03)}
& \ {\scriptsize(6.56e-05)} \ \ \ {\scriptsize(8.15e-02)}  \\
[0.5mm]         \Xhline{0.8\arrayrulewidth}
50,000
& 7.63e-01 \ \ \ 	9.93e-01
& 8.01e-01 \ \ \ 	2.03e-01
& 9.30e-01 \ \ \ 	1.57e-02
& 1.00e+00 \ \ \ 	3.04e-01  \\ [-1.2ex]
&  \ {\scriptsize(4.10e-01)} \ \ \  {\scriptsize(6.07e-03)}
&  \ {\scriptsize(1.19e-01)} \ \ \  {\scriptsize(1.43e-01)}
& \ {\scriptsize(1.58e-02)} \ \ \  {\scriptsize(1.08e-03)}
& \ {\scriptsize(4.92e-03)} \ \ \  {\scriptsize(1.61e-01)}  \\
[0.5mm]         \Xhline{0.8\arrayrulewidth}
100,000
& 6.72e-01 \ \ \	8.10e-01
& 6.21e-01 \ \ \	3.50e-02
& 9.24e-01 \ \ \	1.24e-02
& 9.93e-01 \ \ \	1.67e-01  \\ [-1.2ex]
&  \ {\scriptsize(2.31e-01)} \ \ \  {\scriptsize(2.12e-01)}
&  \ {\scriptsize(1.30e-01)} \ \ \  {\scriptsize(1.95e-02)}
& \ {\scriptsize(4.00e-02)} \ \ \ {\scriptsize(1.67e-03)}
& \ {\scriptsize(1.67e-02)} \ \ \  {\scriptsize(1.43e-01)}  \\
[0.5mm]     \Xhline{3\arrayrulewidth}
& $\sin(x),$ \ $b=8\pi$
& $\sin(x),$ \ $b=16\pi$
& $\sin(x),$ \ $b=32\pi$
& $\sin(x),$ \ $b=64\pi$ \\  [-0.8ex]
Max.\ epoch & $\epsilon_u(\widetilde{X}), \qquad \ \ \epsilon_f(X)$ & $\epsilon_u(\widetilde{X}), \qquad \ \  \epsilon_f(X)$ &
 $\epsilon_u(\widetilde{X}), \qquad \ \  \epsilon_f(X)$ & $\epsilon_u(\widetilde{X}), \qquad \ \ \epsilon_f(X)$ \\[1mm]        \Xhline{0.8\arrayrulewidth}
20,000
& 8.00e-01 \ \ \	1.03e-02
& 7.34e-01 \ \ \	7.08e-03
& 8.05e-01 \ \ \	7.37e-03
& 9.69e-01 \ \ \	8.69e-03  \\ [-1.2ex]
&  \ {\scriptsize(7.03e-02)} \ \ \  {\scriptsize(1.93e-03)}
&  \ {\scriptsize(1.18e-01)} \ \ \  {\scriptsize(1.07e-03)}
& \ {\scriptsize(7.40e-02)} \ \ \ {\scriptsize(6.04e-04)}
& \ {\scriptsize(5.11e-02)} \ \ \  {\scriptsize(4.37e-04)}  \\
[0.5mm]         \Xhline{0.8\arrayrulewidth}
50,000
& 5.45e-01 \ \ \	6.32e-03
& 4.50e-01 \ \ \	4.33e-03
& 5.40e-01 \ \ \	5.02e-03
& 9.15e-01 \ \ \	8.07e-03  \\ [-1.2ex]
&  \ {\scriptsize(9.67e-02)} \ \ \ {\scriptsize(1.50e-03)}
&  \ {\scriptsize(1.49e-01)} \ \ \ {\scriptsize(1.26e-03)}
& \ {\scriptsize(1.13e-01)} \ \ \ {\scriptsize(8.65e-04)}
& \ {\scriptsize(5.27e-02)} \ \ \ {\scriptsize(3.93e-04)}  \\
[0.5mm]         \Xhline{0.8\arrayrulewidth}
100,000
& 2.94e-01 \ \ \	3.73e-03
& 2.02e-01 \ \ \	2.19e-03
& 2.91e-01 \ \ \	2.75e-03
& 7.88e-01 \ \ \	6.95e-03  \\ [-1.2ex]
&  \ {\scriptsize(1.17e-01)} \ \ \ {\scriptsize(1.36e-03)}
&  \ {\scriptsize(1.09e-01)} \ \ \  {\scriptsize(8.24e-04)}
& \ {\scriptsize(8.68e-02)} \ \ \ {\scriptsize(7.47e-04)}
& \ {\scriptsize(9.84e-02)} \ \ \ {\scriptsize(7.52e-04)}  \\
[0.1mm]     \Xhline{3\arrayrulewidth}
    \end{tabular}
}

\vskip-.7truecm
\end{center}
\end{table}

\subsubsection{The residual correction step}
\begin{figure}[htb!]
 \centering
  \includegraphics[width=8cm]{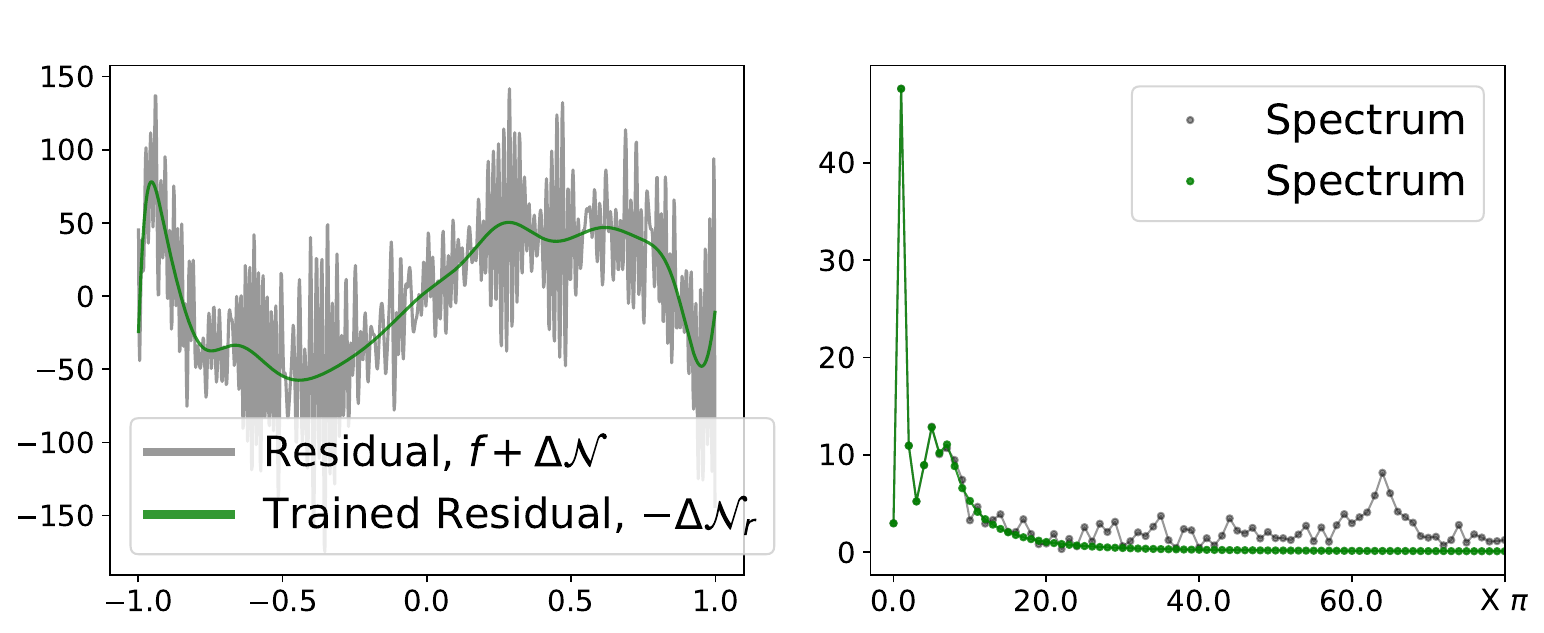}
  \includegraphics[width=7.7cm]{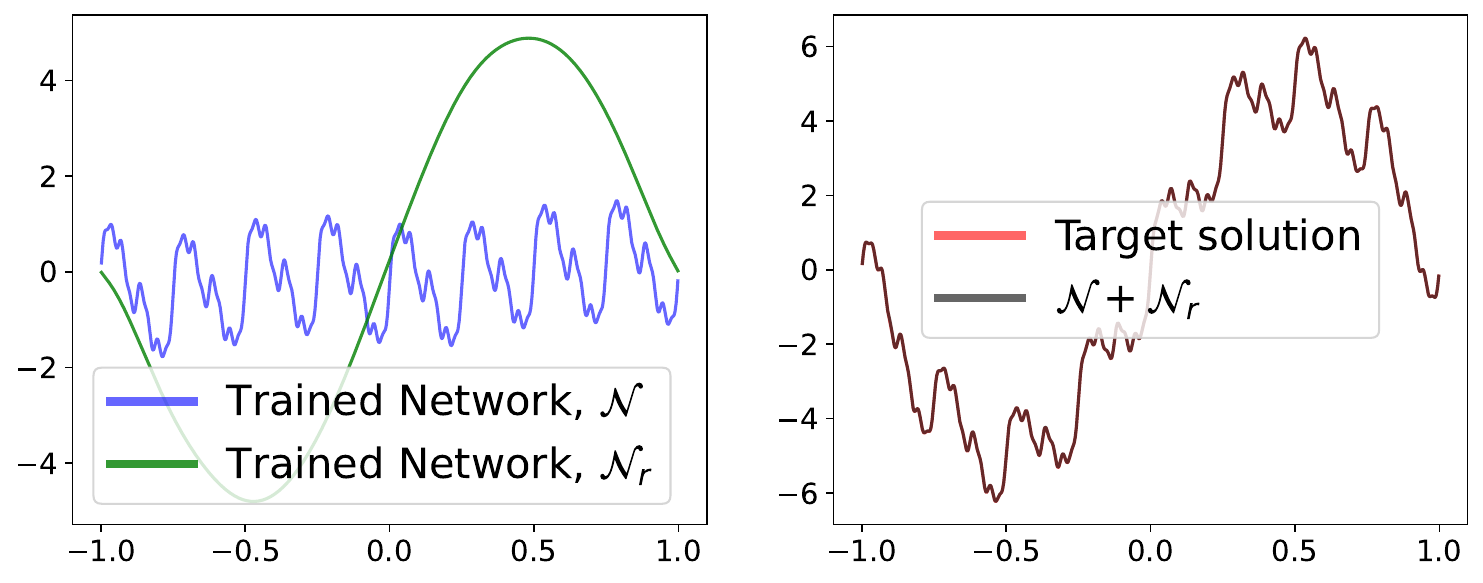}
  \caption{Results of the residual training step. Left two panels: plots of the residual $f+\triangle \mathcal{N}(x;\theta)$ and $ -\triangle \mathcal{N}_r(x;\theta_r)$ and their corresponding frequency spectra.
  Right two panels: plots of the trained network $\mathcal{N}(x;\theta)$ and $\mathcal{N}_r (x;\theta_r)$ and plots of the target solution and the resulting neural network solution
  $\mathcal{N}(x;\theta)+\mathcal{N}_r(x;\theta_r)$.
  The $\sin(5\pi)$ activation function is used in $\mathcal{N}(x;\theta)$.}
  \label{sec2_fig7}
 \end{figure}

We now confirm the effectiveness of the residual correction step at capturing
the remaining low-frequency parts of the error.
As a showcase, we apply it to the trained solution $\mathcal{N}(x;\theta)$ from the previous subsection using the $\sin(5x)$ activation function.
In Fig.~\ref{sec2_fig7}, we plot the residual errors and their frequency spectra for the trained neural network solution $\mathcal{N}(x;\theta)$.
As can be seen in the spectra, the high-frequency parts in the right-hand side function $f(x)$ are resolved well enough.
However, the low-frequency parts of $f(x)$ remain as dominant errors.
To resolve such smooth low-frequency components, we apply the residual correction method with $\mathcal{N}(x;\theta)$ instead of $\mathcal{N}(x;\theta_s)$ in \eqref{pb:residual-correction}.
We perform only 10,000 training epochs to obtain $\mathcal{N}_r(x; \theta_r)$.
Such settings are sufficient for learning smooth and low-frequency feature solutions.

In Fig.~\ref{sec2_fig7}, we also plot the trained residual $-\triangle \mathcal{N}_r(x;\theta_r)$,
its spectrum, and the trained solution $\mathcal{N}_r(x;\theta_r)$.
As can be seen in the spectra, the dominant low-frequency components of the residual error, $f+\triangle\mathcal{N}$,
are well approximated by the trained residual $-\triangle \mathcal{N}_r$.
The solution plots $\mathcal{N}$ and $\mathcal{N}+\mathcal{N}_r$, before and after the residual correction step,
also show that the post-processing step works effectively to capture the smooth part of the error that cannot be resolved by the first-stage neural network approximation function $\mathcal{N}(x;\theta)$ with the $\sin(5x)$
activation function.

We apply the residual training step for the neural network solutions obtained in the previous subsection
and report the final error results in Table~\ref{sec2_tb2}.
The error values are reported for the resulting trained neural network approximation,
$\mathcal{N}(x;\theta)+\mathcal{N}_r(x;\theta_r)$ for the $\sin(ax)$-type activation functions
and $\mathcal{N}(x;\theta_s)+\mathcal{N}_r(x;\theta_r)$ for the scaled Poisson problem with the scaling factor
$b=8\pi,\; 16\pi,\;32\pi,\; 64\pi$.
For the first-step neural network solution, i.e., $\mathcal{N}$, we set the maximum number of epochs to
20,000, 50,000, or 100,000, and the residual error equation is then formed by using the trained neural network solution $\mathcal{N}$.
As seen in the error values in Table~\ref{sec2_tb2}, by finding the first-step neural network solution $\mathcal{N}(x;\theta)$ more accurately with a sufficiently large number of training epochs, low-frequency components become more dominant in the residual error, and the resulting error in $\mathcal{N}+\mathcal{N}_r$ after the residual training step becomes smaller.
In addition, the domain scaling approach gives more accurate approximations than simply using
the $\sin(ax)$-type activation functions. On the other hand, with an excessively large scaling factor $b$, the domain scaling approach may not capture the high-frequency parts of the model solution accurately enough and may produce larger errors than those with smaller scaling factors $b$, as seen in the case with $b=64\pi$.

Based on these experimental results for the one-dimensional model problem,
we can conclude that the domain scaling approach works more effectively
for capturing the high-frequency parts of the model solution than simply using
the $\sin(ax)$-type activation functions, and that 
the remaining low-frequency part of the solution can be well approximated
by the post-processing residual correction step.
In addition, when training the neural network solution for the scaled problem,
it is beneficial to use a larger network size and more training epochs
so that the low-frequency components become more dominant in the residual error.
For training such a dominant low-frequency residual error, we can employ
a smaller network size and fewer training epochs.

\begin{table}[ht!]
\begin{center}
  \caption{
   Error values $\epsilon_u^r(\widetilde{X})$ and $\epsilon_f^r(X)$ for the residual training step applied to the trained solution $\mathcal{N}$ in Table~\ref{sec2_tb1}.
   The residual error equation is formed for the trained solution $\mathcal{N}$ at the maximum epoch.
   The errors are reported for the resulting trained solution, $\mathcal{N}+\mathcal{N}_r$.
   The mean values and standard deviations are listed for six different training sample sets $X$.
   The results are extracted from the spot with the smallest loss value during the training epochs.}\label{sec2_tb2}
{\footnotesize \renewcommand{\arraystretch}{1.2}
    \begin{tabular}{cccccc}
        \Xhline{3\arrayrulewidth}
& $\sin(x)$ & $\sin(2x)$ & $\sin(5x)$ & $\sin(10x)$ \\  [-0.8ex]
Max.\ epoch for $\mathcal{N}$ & $\epsilon_u^r(\widetilde{X}), \qquad \ \ \epsilon_f^r(X)$ & $\epsilon_u^r(\widetilde{X}), \qquad \ \ \epsilon_f^r(X)$ &
 $\epsilon_u^r(\widetilde{X}), \qquad \ \ \epsilon_f^r(X)$ & $\epsilon_u^r(\widetilde{X}), \qquad \ \ \epsilon_f^r(X)$ \\[1mm]         \Xhline{0.8\arrayrulewidth}
20,000
& 1.03e+00 \ \ \	1.00e+00
& 8.74e-01 \ \ \	6.58e-01
& 5.40e-02 \ \ \	2.15e-02
& 1.07e+00 \ \ \	5.08e-01  \\ [-1.2ex]
&  \ {\scriptsize(3.66e-02)} \ \ \ {\scriptsize(1.22e-05)}
&  \ {\scriptsize(3.03e-01)} \ \ \ {\scriptsize(1.36e-01)}
& \ {\scriptsize(3.58e-02)} \ \ \ {\scriptsize(6.10e-03)}
& \ {\scriptsize(1.83e-01)} \ \ \ {\scriptsize(8.23e-02)}  \\
[0.5mm]         \Xhline{0.8\arrayrulewidth}
50,000
& 7.64e-01 \ \ \	9.93e-01
& 4.40e-01 \ \ \	2.02e-01
& 3.04e-02 \ \ \	1.05e-02
& 9.95e-01 \ \ \	2.98e-01  \\ [-1.2ex]
&  \ {\scriptsize(4.10e-01)} \ \ \  {\scriptsize(6.07e-03)}
&  \ {\scriptsize(4.01e-01)} \ \ \  {\scriptsize(1.43e-01)}
& \ {\scriptsize(2.84e-02)} \ \ \  {\scriptsize(1.58e-03)}
& \ {\scriptsize(3.46e-01)} \ \ \ {\scriptsize(1.68e-01)}  \\
[0.5mm]         \Xhline{0.8\arrayrulewidth}
100,000
& 6.74e-01 \ \ \	8.10e-01
& 4.22e-02 \ \ \	3.29e-02
& 1.64e-02 \ \ \	6.87e-03
& 1.57e+00 \ \ \	1.50e-01  \\ [-1.2ex]
&  \ {\scriptsize(2.36e-01)} \ \ \  {\scriptsize(2.12e-01)}
&  \ {\scriptsize(2.71e-02)} \ \ \  {\scriptsize(2.00e-02)}
& \ {\scriptsize(1.13e-02)} \ \ \ {\scriptsize(1.58e-03)}
& \ {\scriptsize(9.06e-01)} \ \ \  {\scriptsize(1.44e-01)}  \\
[0.5mm]     \Xhline{3\arrayrulewidth}
& $\sin(x),$ \ $b=8\pi$
& $\sin(x),$ \ $b=16\pi$
& $\sin(x),$ \ $b=32\pi$
& $\sin(x),$ \ $b=64\pi$ \\  [-0.8ex]
Max.\ epoch for $\mathcal{N}$ & $\epsilon_u^r(\widetilde{X}), \qquad \ \ \epsilon_f^r(X)$ & $\epsilon_u^r(\widetilde{X}), \qquad \ \ \epsilon_f^r(X)$ &
 $\epsilon_u^r(\widetilde{X}), \qquad \ \ \epsilon_f^r(X)$ & $\epsilon_u^r(\widetilde{X}), \qquad \ \ \epsilon_f^r(X)$ \\[1mm]         \Xhline{0.8\arrayrulewidth}
20,000
& 1.36e-02 \ \ \	7.18e-03
& 8.44e-03 \ \ \	2.96e-03
& 8.36e-03 \ \ \	2.23e-03
& 1.58e-02 \ \ \ 	2.04e-03  \\ [-1.2ex]
&  \ {\scriptsize(9.87e-03)} \ \ \  {\scriptsize(2.77e-03)}
&  \ {\scriptsize(9.14e-03)} \ \ \  {\scriptsize(3.34e-04)}
& \ {\scriptsize(4.94e-03)} \ \ \  {\scriptsize(2.66e-04)}
& \ {\scriptsize(1.05e-02)} \ \ \  {\scriptsize(3.16e-04)}  \\
[0.5mm]         \Xhline{0.8\arrayrulewidth}
50,000
& 1.35e-02 \ \ \	4.04e-03
& 8.85e-03 \ \ \	1.75e-03
& 7.05e-03 \ \ \	1.68e-03
& 1.52e-02 \ \ \	1.60e-03  \\ [-1.2ex]
&  \ {\scriptsize(1.12e-02)} \ \ \  {\scriptsize(1.55e-03)}
&  \ {\scriptsize(1.01e-02)} \ \ \  {\scriptsize(2.77e-04)}
& \ {\scriptsize(4.12e-03)} \ \ \  {\scriptsize(3.24e-04)}
& \ {\scriptsize(8.86e-03)} \ \ \  {\scriptsize(2.99e-04)}  \\
[0.5mm]         \Xhline{0.8\arrayrulewidth}
100,000
& 8.21e-03 \ \ \	2.59e-03
& 5.62e-03 \ \ \	1.15e-03
& 4.99e-03 \ \ \	1.03e-03
& 1.44e-02 \ \ \	1.41e-03  \\ [-1.2ex]
&  \ {\scriptsize(7.57e-03)} \ \ \  {\scriptsize(1.11e-03)}
&  \ {\scriptsize(9.25e-03)} \ \ \  {\scriptsize(2.18e-04)}
& \ {\scriptsize(2.99e-03)} \ \ \  {\scriptsize(1.84e-04)}
& \ {\scriptsize(8.76e-03)} \ \ \  {\scriptsize(1.97e-04)}  \\
[0.1mm]     \Xhline{3\arrayrulewidth}
    \end{tabular}
}

\vskip-.7truecm
\end{center}
\end{table}

\subsection{Two-dimensional Poisson problem}
\begin{figure}[htb!]
 \centering
  \includegraphics[height=3.4cm]{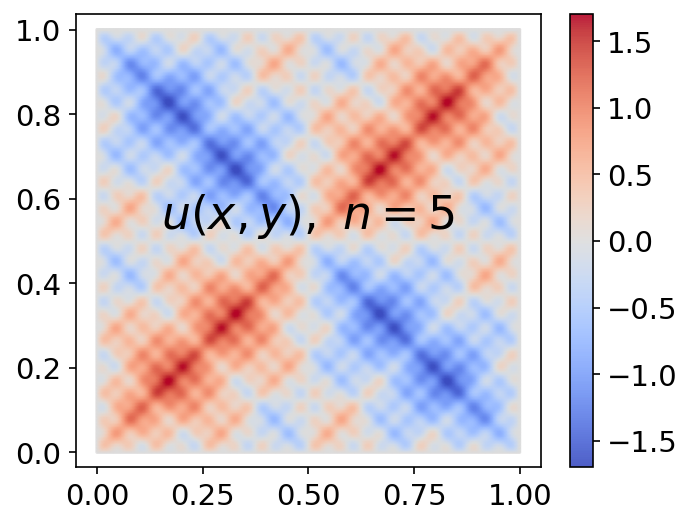}
  \includegraphics[height=3.4cm]{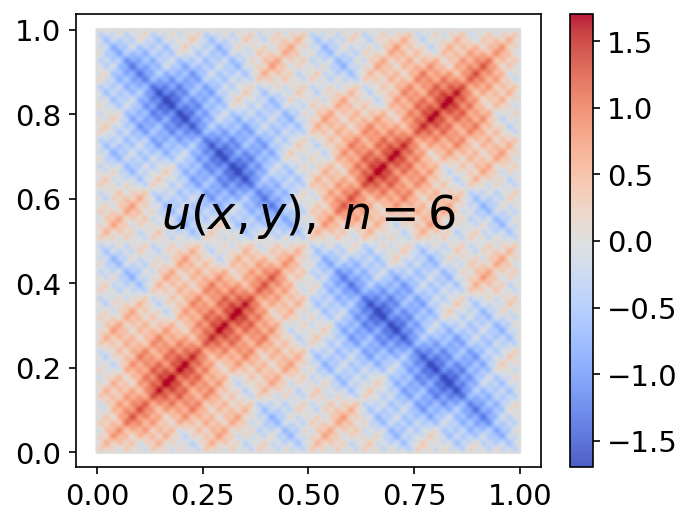}
  \caption{The target solution $u(x,y)$ with $n=5$ (left) and $n=6$ (right).}
  \label{sec22_fig1}
\end{figure}

We next consider a two-dimensional Poisson problem with a multi-frequency solution
\begin{equation}\label{sol_2d}
u(x,y) = \sum_{i=1}^n \frac{1}{i} \sin(2^i \pi x) \sin(2^i \pi y)  
\end{equation}
in $\Omega =(0, 1) \times (0, 1)$ for $n=5$ and $6$.
We present plots of the model solution for both $n=5$ and $6$ in Fig.~\ref{sec22_fig1}.
The solutions include highly multiscale features, so it is challenging to find a good approximate solution using the neural network approximation in this case.
We can apply our two-step strategy to such multiscale feature problems.

\subsubsection{The case $n=5$}
As mentioned above, we set up the network $\mathcal{N}$ as a fully connected neural network with four hidden layers and with 70 neurons per each hidden layer.
To form a training data set $X$, we choose 160,000 points from the domain interior and 4,000 points from the boundary.
For the residual correction network $\mathcal{N}_r$, we set it up as a fully connected neural network with four hidden layers and 50 neurons in each hidden layer.
We perform 10,000 epochs to train the network parameters in $\mathcal{N}_r$.

\begin{figure}[htb!]
 \centering
  \includegraphics[height=2.9cm]{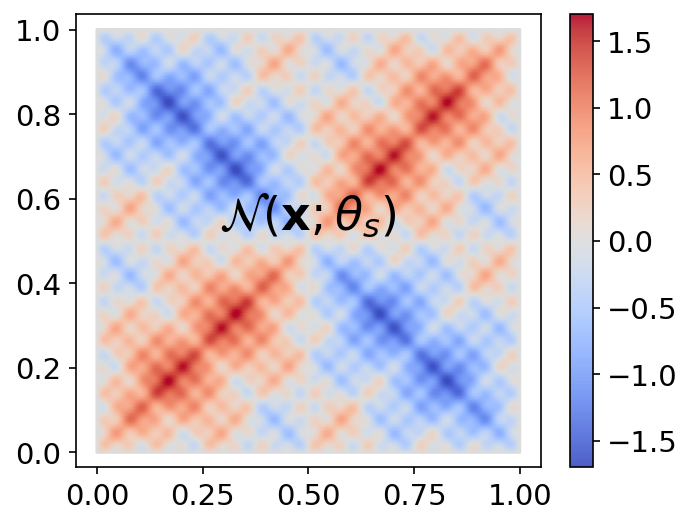}
  \includegraphics[height=2.9cm]{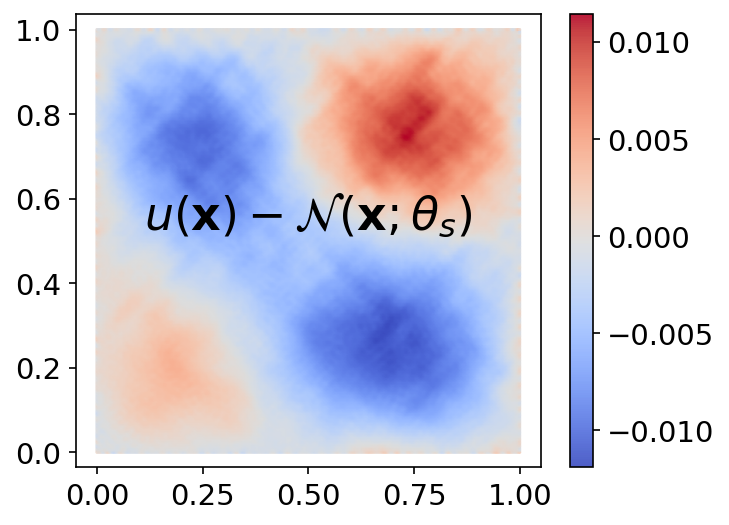}
  \includegraphics[height=2.9cm]{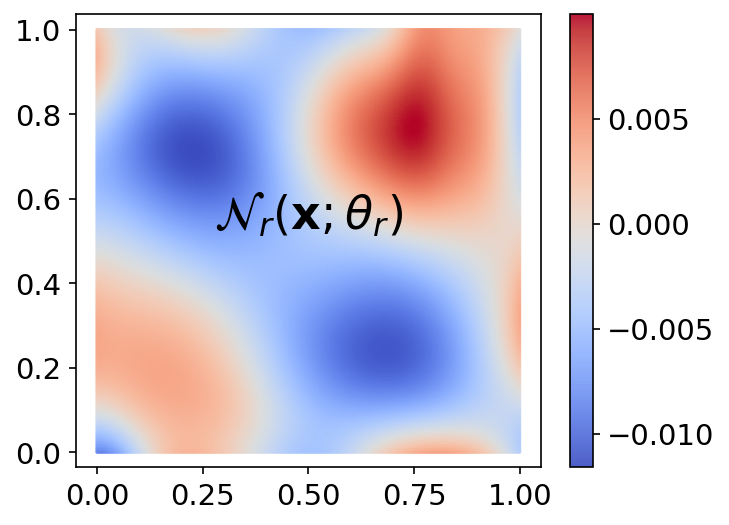}
  \includegraphics[height=2.9cm]{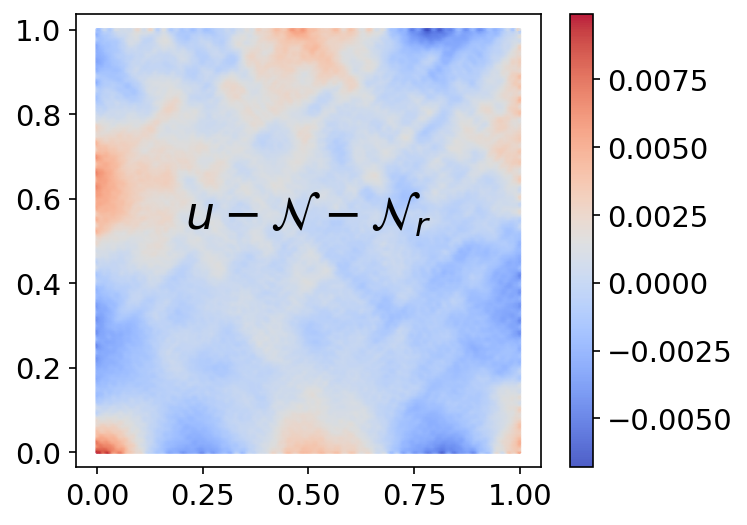}
  \caption{Multi-frequency model solution to \eqref{sol_2d} with $n=5$: plots of the trained solution $\mathcal{N}$ with $b=32\pi$, the error $u - \mathcal{N}$, the residual trained solution $\mathcal{N}_r$, and the error in the resulting solution, $u-\mathcal{N}-\mathcal{N}_r$.}
  \label{sec22_fig2}
\end{figure}

To demonstrate the effectiveness of our strategy,
we plot the trained solution $\mathcal{N}(\mathbf{x};\theta_s)$ from the domain scaling procedure with $b=32\pi$,
the residual trained solution $\mathcal{N}_r(\mathbf{x};\theta_r)$, and the errors in the resulting solution, $\mathcal{N}+\mathcal{N}_r$ in Fig.~\ref{sec22_fig2}.
We can see the residual trained solution $\mathcal{N}_r ( \mathbf{x}; \theta_r)$ approximates the target residual $u(\mathbf{x}) - \mathcal{N}(\mathbf{x};\theta_s)$ well. The target solution $u(\mathbf{x})$ is thus accurately resolved by the combination $\mathcal{N}(\mathbf{x};\theta_s) + \mathcal{N}_r (\mathbf{x}; \theta_r)$.

In Table~\ref{sec22_tb1}, we list the error values for various scaling factors $b=8\pi, 16\pi$, and $32\pi$ and with the maximum number of training epochs set to 20,000, 50,000, and 100,000.
From the $\epsilon_f$ values, we can see that the domain scaling method effectively approximates the high-frequency components, while values of $\epsilon_u$ larger than $\epsilon_f$ present the need
for the additional residual correction step
to improve the errors in the low-frequency components.
The residual correction step approximates the low-frequency components that could not be resolved in the scaled problem, leading to smaller error values for ${\epsilon}_u^r$ compared to ${\epsilon}_u$.
As already observed in the one-dimensional example,
by training the network solution $\mathcal{N}$ for the scaled problem with a large enough number of training epochs, around 100,000 in this case,
we can obtain a residual error equation that gives a more accurate solution, $\mathcal{N}+\mathcal{N}_r$,
after the residual correction step.
\begin{table}[ht!]
\begin{center}
  \caption{
Error values for multi-frequency model problem \eqref{sol_2d} with $n=5$: $\epsilon_u(\widetilde{X})$ and $\epsilon_f(X)$ for the trained solution $\mathcal{N}$ of the scaled problem and $\epsilon_u^r(\widetilde{X})$ and $\epsilon_f^r(X)$
 for the resulting solution $\mathcal{N}+\mathcal{N}_r$ after the residual correction step.
 The mean values and standard deviations are listed for six different training sample sets $X$. The results are extracted from the spot with the smallest loss value during the training epochs.}\label{sec22_tb1}
{\footnotesize \renewcommand{\arraystretch}{1.2}
    \begin{tabular}{cccccc}
        \Xhline{3\arrayrulewidth}
$n=5$
& $\sin(x),$ \ $b=8\pi$
& $\sin(x),$ \ $b=16\pi$
& $\sin(x),$ \ $b=32\pi$ \\  [-0.8ex]
Max.\ epoch for $\mathcal{N}$  & $\epsilon_u(\widetilde{X}), \qquad \ \  \epsilon_f(X)$ &
 $\epsilon_u(\widetilde{X}), \qquad \ \  \epsilon_f(X)$ & $\epsilon_u(\widetilde{X}), \qquad \ \  \epsilon_f(X)$ \\[1mm]         \Xhline{0.8\arrayrulewidth}
20,000
& 3.34e-02	\ \ \ 2.11e-02
& 3.08e-02 \ \ \	1.14e-02
& 2.92e-02 \ \ \	8.32e-03  \\ [-1.2ex]
&  \ {\scriptsize(1.12e-02)} \ \ \  {\scriptsize(2.86e-03)}
& \ {\scriptsize(2.63e-02)} \ \ \  {\scriptsize(1.02e-03)}
& \ {\scriptsize(1.15e-02)} \ \ \  {\scriptsize(3.18e-04)}  \\
[0.5mm]         \Xhline{0.8\arrayrulewidth}
50,000
& 1.53e-02 \ \ \ 	1.18e-02
& 1.37e-02 \ \ \ 	6.73e-03
& 1.58e-02 \ \ \ 	5.18e-03  \\ [-1.2ex]
&  \ {\scriptsize(1.50e-02)} \ \ \  {\scriptsize(1.54e-03)}
& \ {\scriptsize(5.85e-03)} \ \ \  {\scriptsize(6.68e-04)}
& \ {\scriptsize(4.67e-03)} \ \ \  {\scriptsize(2.32e-04)}  \\
[0.5mm]         \Xhline{0.8\arrayrulewidth}
100,000
& 9.85e-03 \ \ \ 	7.35e-03
& 6.85e-03 \ \ \ 	4.57e-03
& 1.04e-02 \ \ \ 	3.81e-03  \\ [-1.2ex]
&  \ {\scriptsize(8.06e-03)} \ \ \  {\scriptsize(9.31e-04)}
& \ {\scriptsize(2.91e-03)} \ \ \  {\scriptsize(5.79e-04)}
& \ {\scriptsize(1.79e-03)} \ \ \  {\scriptsize(2.27e-04)}  \\
[0.1mm]     \Xhline{3\arrayrulewidth}
Max.\ epoch for $\mathcal{N}$  & $\epsilon_u^r(\widetilde{X}), \qquad \ \  \epsilon_f^r(X)$ &
 $\epsilon_u^r(\widetilde{X}), \qquad \ \  \epsilon_f^r(X)$ & $\epsilon_u^r(\widetilde{X}), \qquad \ \  \epsilon_f^r(X)$ \\[1mm]         \Xhline{0.8\arrayrulewidth}
20,000
& 1.02e-02	 \ \ \ 2.11e-02
& 6.55e-03	 \ \ \ 1.14e-02
& 5.73e-03	 \ \ \ 8.31e-03  \\ [-1.2ex]
&  \ {\scriptsize(4.80e-03)} \ \ \  {\scriptsize(2.86e-03)}
& \ {\scriptsize(2.69e-03)} \ \ \  {\scriptsize(1.02e-03)}
& \ {\scriptsize(3.32e-03)} \ \ \  {\scriptsize(3.22e-04)}  \\
[0.5mm]         \Xhline{0.8\arrayrulewidth}
50,000
& 6.25e-03	 \ \ \ 1.18e-02
& 3.27e-03	 \ \ \ 6.73e-03
& 3.92e-03	 \ \ \ 5.17e-03  \\ [-1.2ex]
&  \ {\scriptsize(2.86e-03)} \ \ \  {\scriptsize(1.54e-03)}
& \ {\scriptsize(1.33e-03)} \ \ \  {\scriptsize(6.70e-04)}
& \ {\scriptsize(9.09e-04)} \ \ \  {\scriptsize(2.34e-04)}  \\
[0.5mm]         \Xhline{0.8\arrayrulewidth}
100,000
& 3.50e-03	 \ \ \ 7.35e-03
& 2.41e-03	 \ \ \ 4.56e-03
& 3.47e-03	 \ \ \ 3.80e-03  \\ [-1.2ex]
&  \ {\scriptsize(1.92e-03)} \ \ \  {\scriptsize(9.31e-04)}
& \ {\scriptsize(8.65e-04)} \ \ \  {\scriptsize(5.78e-04)}
& \ {\scriptsize(1.25e-03)} \ \ \  {\scriptsize(2.26e-04)}  \\
[0.1mm]     \Xhline{3\arrayrulewidth}
    \end{tabular}
}

\vskip-.7truecm
\end{center}
\end{table}

\subsubsection{The case $n=6$}
We finally consider a more challenging multi-frequency model solution with $n=6$.
For this case, we set a larger network size and training data set for the scaled problem.
The network function $\mathcal{N}$ is set to a fully connected neural network with four hidden layers, each with 100 neurons. The training data set $X$ consists of 250,000 points chosen from the domain interior and 6,000 points chosen from the boundary.
For the residual training, we set the network function $\mathcal{N}_r$ with four hidden layers and 50 neurons per hidden layer, and train it with 10,000 epochs.

\begin{figure}[htb!]
 \centering
  \includegraphics[height=2.9cm]{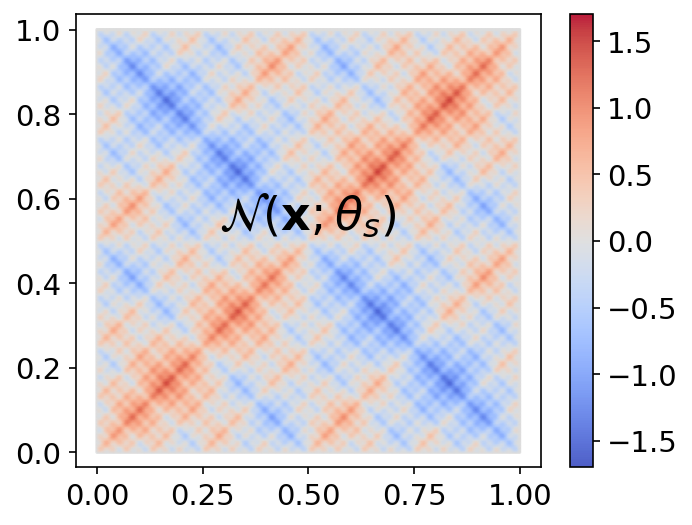}
  \includegraphics[height=2.9cm]{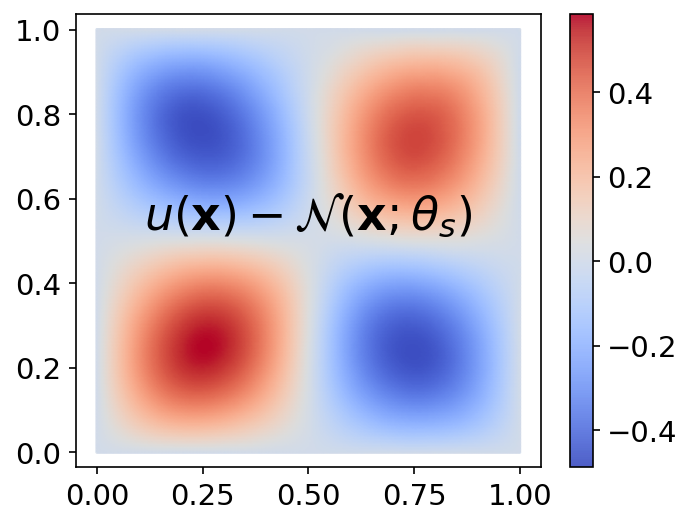}
  \includegraphics[height=2.9cm]{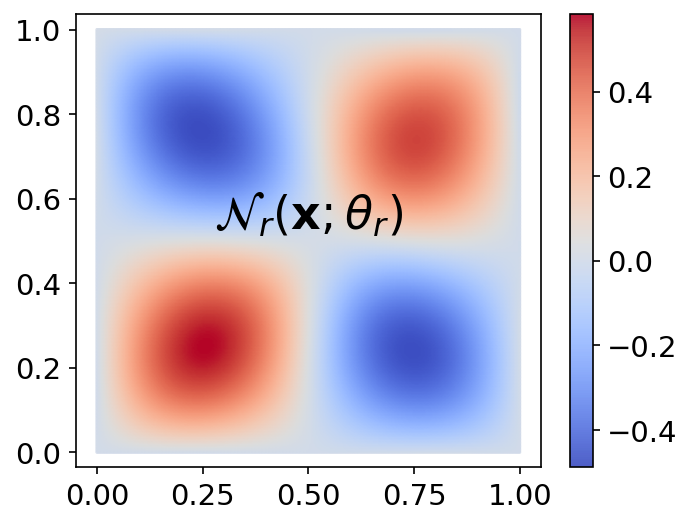}
  \includegraphics[height=2.9cm]{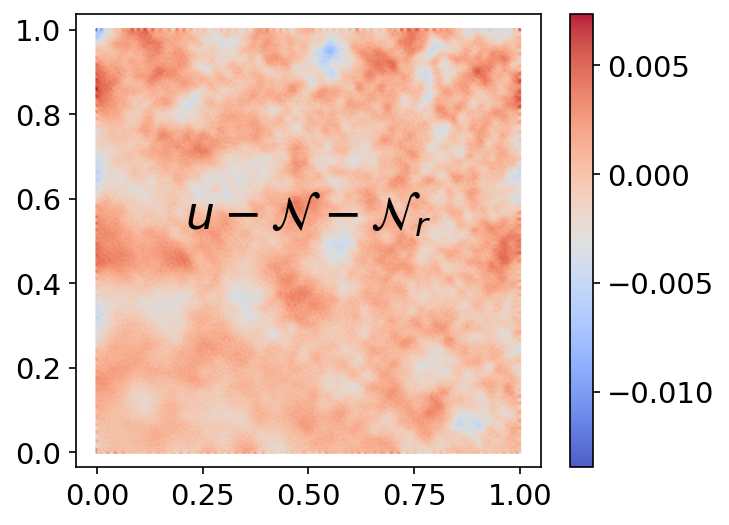}
  \caption{
  Multi-frequency model solution to \eqref{sol_2d} with $n=6$: plots of trained solution $\mathcal{N}$ with $b=64\pi$, the error $u - \mathcal{N}$, the residual trained solution $\mathcal{N}_r$, and the error in the resulting solution, $u-\mathcal{N}-\mathcal{N}_r$.}
  \label{sec22_fig3}
\end{figure}

As a showcase, we present a plot of the trained solution $\mathcal{N}(\mathbf{x};\theta_s)$ with the scaling factor $b=64\pi$ after 100,000 training epochs, as well as plots of the error $u(\mathbf{x})-\mathcal{N}(\mathbf{x};\theta_s)$, the residual trained solution $\mathcal{N}_r(\mathbf{x};\theta_r)$, and the error in the resulting solution, in Fig.~\ref{sec22_fig3}.
The trained solution $\mathcal{N}(\mathbf{x};\theta_s)$ for the scaled problem
captures the high-frequency components of the target solution effectively.
However, it still exhibits a considerable amount of low-frequency error as seen in the error plot.
The outcome of the residual correction process demonstrates that $\mathcal{N}_r ( \mathbf{x}; \theta_r)$ effectively approximates the remaining error $u(\mathbf{x}) - \mathcal{N}(\mathbf{x};\theta_s)$.
Consequently, the target function is well approximated by the combination $\mathcal{N}(\mathbf{x};\theta_s) + \mathcal{N}_r ( \mathbf{x}; \theta_r)$.

We report the error values of our proposed methods for scaling factor $b=16\pi$, $32\pi$, and $64\pi$, and with the maximum number of training epochs set to 20,000, 50,000, and 100,000 for the neural network solution $\mathcal{N}$ of the scaled model problem, in Table~\ref{sec22_tb2}.
The error values are computed for the scaled model solution $\mathcal{N}$ and for the resulting solution $\mathcal{N}+\mathcal{N}_r$ after the residual correction step.
From the $\epsilon_f$ values, we can again see that the domain scaling method effectively approximates the high-frequency components. In addition, the residual training step works successfully captures the remaining
low-frequency parts as seen in the error values $\epsilon_u^r$ and $\epsilon_u$.

\begin{table}[ht!]
\begin{center}
  \caption{
  Error values for multi-frequency model problem \eqref{sol_2d} with $n=6$: $\epsilon_u(\widetilde{X})$ and $\epsilon_f(X)$ for the trained solution $\mathcal{N}$ of the scaled problem and $\epsilon_u^r(\widetilde{X})$ and $\epsilon_f^r(X)$
 for the resulting solution $\mathcal{N}+\mathcal{N}_r$ after the residual correction step.
 The mean values and standard deviations are listed for six different training sample sets $X$. The results are extracted from the spot with the smallest loss value during the training epochs.}\label{sec22_tb2}
{\footnotesize \renewcommand{\arraystretch}{1.2}
    \begin{tabular}{cccccc}
        \Xhline{3\arrayrulewidth}
$n=6$
& $\sin(x),$ \ $b=16\pi$
& $\sin(x),$ \ $b=32\pi$
& $\sin(x),$ \ $b=64\pi$ \\  [-0.8ex]
Max.\ epoch for $\mathcal{N}$  & $\epsilon_u(\widetilde{X}), \qquad \ \  \epsilon_f(X)$ &
 $\epsilon_u(\widetilde{X}), \qquad \ \  \epsilon_f(X)$ & $\epsilon_u(\widetilde{X}), \qquad \ \  \epsilon_f(X)$ \\[1mm]         \Xhline{0.8\arrayrulewidth}
20,000
& 3.82e-01	 \ \ \ 4.14e-02
& 3.15e-01	 \ \ \ 2.13e-02
& 5.92e-01	 \ \ \ 1.77e-02  \\ [-1.2ex]
&  \ {\scriptsize(9.63e-02)} \ \ \  {\scriptsize(1.94e-03)}
& \ {\scriptsize(1.25e-01)} \ \ \  {\scriptsize(1.71e-03)}
& \ {\scriptsize(6.40e-02)} \ \ \  {\scriptsize(7.17e-04)}  \\
[0.5mm]         \Xhline{0.8\arrayrulewidth}
50,000
& 1.72e-01	 \ \ \ 2.23e-02
& 1.68e-01	 \ \ \ 1.29e-02
& 4.79e-01	 \ \ \ 1.27e-02  \\ [-1.2ex]
&  \ {\scriptsize(1.08e-01)} \ \ \  {\scriptsize(9.26e-04)}
& \ {\scriptsize(6.42e-02)} \ \ \  {\scriptsize(6.87e-04)}
& \ {\scriptsize(6.51e-02)} \ \ \  {\scriptsize(3.47e-04)}  \\
[0.5mm]         \Xhline{0.8\arrayrulewidth}
100,000
& 8.66e-02		 \ \ \ 1.46e-02
& 1.26e-01		 \ \ \ 9.42e-03
& 4.04e-01		 \ \ \ 1.02e-02  \\ [-1.2ex]
&  \ {\scriptsize(5.19e-02)} \ \ \  {\scriptsize(8.00e-04)}
& \ {\scriptsize(2.83e-02)} \ \ \  {\scriptsize(4.28e-04)}
& \ {\scriptsize(5.51e-02)} \ \ \  {\scriptsize(4.60e-04)}  \\
[0.1mm]     \Xhline{3\arrayrulewidth}
Max.\ epoch for $\mathcal{N}$  & $\epsilon_u^r(\widetilde{X}), \qquad \ \  \epsilon_f^r(X)$ &
 $\epsilon_u^r(\widetilde{X}), \qquad \ \  \epsilon_f^r(X)$ & $\epsilon_u^r(\widetilde{X}), \qquad \ \  \epsilon_f^r(X)$ \\[1mm]         \Xhline{0.8\arrayrulewidth}
20,000
& 6.00e-02	 \ \ \ 4.14e-02
& 4.32e-02	 \ \ \ 2.12e-02
& 5.22e-02	 \ \ \ 1.73e-02  \\ [-1.2ex]
&  \ {\scriptsize(4.18e-02)} \ \ \  {\scriptsize(1.95e-03)}
& \ {\scriptsize(1.92e-02)} \ \ \  {\scriptsize(1.71e-03)}
& \ {\scriptsize(2.41e-02)} \ \ \  {\scriptsize(8.14e-04)}  \\
[0.5mm]         \Xhline{0.8\arrayrulewidth}
50,000
& 2.95e-02	 \ \ \ 2.23e-02
& 2.56e-02	 \ \ \ 1.29e-02
& 4.58e-02	 \ \ \ 1.24e-02  \\ [-1.2ex]
&  \ {\scriptsize(1.20e-02)} \ \ \  {\scriptsize(1.04e-03)}
& \ {\scriptsize(1.14e-02)} \ \ \  {\scriptsize(6.95e-04)}
& \ {\scriptsize(2.35e-02)} \ \ \  {\scriptsize(3.76e-04)}  \\
[0.5mm]         \Xhline{0.8\arrayrulewidth}
100,000
& 1.67e-02	 \ \ \ 1.46e-02
& 1.93e-02	 \ \ \ 9.39e-03
& 4.28e-02	 \ \ \ 9.93e-03  \\ [-1.2ex]
&  \ {\scriptsize(7.51e-03)} \ \ \  {\scriptsize(7.97e-04)}
& \ {\scriptsize(9.99e-03)} \ \ \  {\scriptsize(4.29e-04)}
& \ {\scriptsize(2.31e-02)} \ \ \  {\scriptsize(4.42e-04)}  \\
[0.1mm]     \Xhline{3\arrayrulewidth}
    \end{tabular}
}

\vskip-.7truecm
\end{center}
\end{table}

\section{Conclusion}
We have proposed effective methods, domain scaling and residual training, for solving multi-frequency PDEs by  the neural network approximation. Our simple two-step approach was numerically confirmed to give accurate and efficient solutions when applied to highly challenging multi-frequency model problems.
By transforming the original problem into a scaled problem, the high-frequency components of the model solution become less oscillatory and can be well trained by neural network approximation. The low-frequency components of the model solution remain as the dominant error and are trained in the second step to give the resultant solution as a sum of the two neural network functions. 
Applications of the proposed methods to more challenging problems, like Helmholtz equations and nonlinear PDEs, will need further investigation into the hyperparameter settings and the residual correction step. 

\bibliographystyle{plain}
\bibliography{pinn_ms}

\begin{thebibliography}{10}

\bibitem{cai2020phase}
Wei Cai, Xiaoguang Li, and Lizuo Liu.
\newblock A phase shift deep neural network for high frequency approximation
  and wave problems.
\newblock {\em SIAM Journal on Scientific Computing}, 42(5):A3285--A3312, 2020.

\bibitem{kingma2014adam}
Diederik~P Kingma and Jimmy Ba.
\newblock Adam: A method for stochastic optimization.
\newblock {\em arXiv preprint arXiv:1412.6980}, 2014.

\bibitem{lagaris1998artificial}
Isaac~E Lagaris, Aristidis Likas, and Dimitrios~I Fotiadis.
\newblock Artificial neural networks for solving ordinary and partial
  differential equations.
\newblock {\em IEEE Transactions on Neural Networks}, 9(5):987--1000, 1998.

\bibitem{li2020multi}
Xi-An Li, Zhi-Qin~John Xu, and Lei Zhang.
\newblock A multi-scale {DNN} algorithm for nonlinear elliptic equations with
  multiple scales.
\newblock {\em Communications in Computational Physics}, 28(5):1886--1906,
  2020.

\bibitem{liu2020multi}
Ziqi Liu, Wei Cai, and Zhi-Qin~John Xu.
\newblock Multi-scale deep neural network {(MscaleDNN)} for solving
  {Poisson-Boltzmann} equation in complex domains.
\newblock {\em Communications in Computational Physics}, 28(5):1970--2001,
  2020.

\bibitem{luo2021frequency}
Tao Luo, Zheng Ma, Zhi-Qin~John Xu, and Yaoyu Zhang.
\newblock Theory of the frequency principle for general deep neural networks.
\newblock {\em CSIAM Transactions on Applied Mathematics}, 2(3):484--507, 2021.

\bibitem{mathews2021uncovering}
Abhilash Mathews, Manaure Francisquez, Jerry~W Hughes, David~R Hatch, Ben Zhu,
  and Barrett~N Rogers.
\newblock Uncovering turbulent plasma dynamics via deep learning from partial
  observations.
\newblock {\em Physical Review E}, 104(2):025205, 2021.

\bibitem{mccormick1987multigrid}
Stephen~F McCormick.
\newblock {\em Multigrid methods}.
\newblock SIAM, 1987.

\bibitem{misra2019mish}
Diganta Misra.
\newblock Mish: A self regularized non-monotonic activation function.
\newblock {\em arXiv preprint arXiv:1908.08681}, 2019.

\bibitem{rahaman2019spectral}
Nasim Rahaman, Aristide Baratin, Devansh Arpit, Felix Draxler, Min Lin, Fred
  Hamprecht, Yoshua Bengio, and Aaron Courville.
\newblock On the spectral bias of neural networks.
\newblock In {\em International Conference on Machine Learning}, pages
  5301--5310. PMLR, 2019.

\bibitem{raissi2019physics}
Maziar Raissi, Paris Perdikaris, and George~E Karniadakis.
\newblock Physics-informed neural networks: A deep learning framework for
  solving forward and inverse problems involving nonlinear partial differential
  equations.
\newblock {\em Journal of Computational physics}, 378:686--707, 2019.

\bibitem{sitzmann2020implicit}
Vincent Sitzmann, Julien Martel, Alexander Bergman, David Lindell, and Gordon
  Wetzstein.
\newblock Implicit neural representations with periodic activation functions.
\newblock {\em Advances in Neural Information Processing Systems},
  33:7462--7473, 2020.

\bibitem{wang2020multi}
Bo~Wang, Wenzhong Zhang, and Wei Cai.
\newblock Multi-scale deep neural network {(MscaleDNN)} methods for oscillatory
  {Stokes} flows in complex domains.
\newblock {\em Communications in Computational Physics}, 28(5):2139--2157,
  2020.

\bibitem{xu2019frequency}
Zhi-Qin~John Xu, Yaoyu Zhang, Tao Luo, Yanyang Xiao, and Zheng Ma.
\newblock Frequency principle: Fourier analysis sheds light on deep neural
  networks.
\newblock {\em Communications in Computational Physics}, 28(5):1746--1767,
  2020.

\bibitem{xu2019training}
Zhi-Qin~John Xu, Yaoyu Zhang, and Yanyang Xiao.
\newblock Training behavior of deep neural network in frequency domain.
\newblock In {\em Neural Information Processing: 26th International Conference,
  ICONIP 2019, Sydney, NSW, Australia, December 12--15, 2019, Proceedings, Part
  I 26}, pages 264--274. Springer, 2019.

\end{thebibliography}

\end{document}